\documentclass[12pt]{article}
\usepackage{latexsym, amssymb}
\DeclareMathAlphabet\gothic{U}{euf}{m}{n}

\makeatletter
\def\eqnarray{\stepcounter{equation}\let\@currentlabel=\theequation
\global\@eqnswtrue
\tabskip\@centering\let\\=\@eqncr
$$\halign to \displaywidth\bgroup\hfil\global\@eqcnt\z@
  $\displaystyle\tabskip\z@{##}$&\global\@eqcnt\@ne
  \hfil$\displaystyle{{}##{}}$\hfil
  &\global\@eqcnt\tw@ $\displaystyle{##}$\hfil
  \tabskip\@centering&\llap{##}\tabskip\z@\cr}

\def\endeqnarray{\@@eqncr\egroup
      \global\advance\c@equation\m@ne$$\global\@ignoretrue}

\def\@yeqncr{\@ifnextchar [{\@xeqncr}{\@xeqncr[5pt]}}
\makeatother

\begin{document}
\bibliographystyle{tom}

\newtheorem{lemma}{Lemma}[section]
\newtheorem{thm}[lemma]{Theorem}
\newtheorem{cor}[lemma]{Corollary}
\newtheorem{voorb}[lemma]{Example}
\newtheorem{rem}[lemma]{Remark}
\newtheorem{prop}[lemma]{Proposition}
\newtheorem{stat}[lemma]{{\hspace{-5pt}}}
\newtheorem{obs}[lemma]{Observation}
\newtheorem{defin}[lemma]{Definition}

\newenvironment{remarkn}{\begin{rem} \rm}{\end{rem}}
\newenvironment{exam}{\begin{voorb} \rm}{\end{voorb}}
\newenvironment{defn}{\begin{defin} \rm}{\end{defin}}
\newenvironment{obsn}{\begin{obs} \rm}{\end{obs}}

\newcommand{\gota}{\gothic{a}}
\newcommand{\gotb}{\gothic{b}}
\newcommand{\gotc}{\gothic{c}}
\newcommand{\gote}{\gothic{e}}
\newcommand{\gotf}{\gothic{f}}
\newcommand{\gotg}{\gothic{g}}
\newcommand{\gothh}{\gothic{h}}
\newcommand{\gotk}{\gothic{k}}
\newcommand{\gotm}{\gothic{m}}
\newcommand{\gotn}{\gothic{n}}
\newcommand{\gotp}{\gothic{p}}
\newcommand{\gotq}{\gothic{q}}
\newcommand{\gotr}{\gothic{r}}
\newcommand{\gots}{\gothic{s}}
\newcommand{\gotu}{\gothic{u}}
\newcommand{\gotv}{\gothic{v}}
\newcommand{\gotw}{\gothic{w}}
\newcommand{\gotz}{\gothic{z}}
\newcommand{\gotA}{\gothic{A}}
\newcommand{\gotB}{\gothic{B}}
\newcommand{\gotG}{\gothic{G}}
\newcommand{\gotL}{\gothic{L}}
\newcommand{\gotS}{\gothic{S}}
\newcommand{\gotT}{\gothic{T}}

\newcounter{teller}
\renewcommand{\theteller}{\Roman{teller}}
\newenvironment{tabel}{\begin{list}%
{\rm \bf \Roman{teller}.\hfill}{\usecounter{teller} \leftmargin=1.1cm
\labelwidth=1.1cm \labelsep=0cm \parsep=0cm}
                      }{\end{list}}

\newcounter{tellerr}
\renewcommand{\thetellerr}{(\roman{tellerr})}
\newenvironment{subtabel}{\begin{list}%
{\rm  (\roman{tellerr})\hfill}{\usecounter{tellerr} \leftmargin=1.1cm
\labelwidth=1.1cm \labelsep=0cm \parsep=0cm}
                         }{\end{list}}
\newenvironment{ssubtabel}{\begin{list}%
{\rm  (\roman{tellerr})\hfill}{\usecounter{tellerr} \leftmargin=1.1cm
\labelwidth=1.1cm \labelsep=0cm \parsep=0cm \topsep=1.5mm}
                         }{\end{list}}

\newenvironment{tabelpairs}{\begin{list}%
{{\rm \bf \Roman{teller}.}\ {\rm( {\bf \Roman{teller}$^\prime.\,$})}\hfill}{\usecounter{teller} \leftmargin=1.9cm
\labelwidth=1.9cm \labelsep=0cm \parsep=0cm}
                     }{\end{list}}

\newcommand{\Ni}{{\bf N}}
\newcommand{\Ri}{{\bf R}}
\newcommand{\Ci}{{\bf C}}
\newcommand{\Ti}{{\bf T}}
\newcommand{\Zi}{{\bf Z}}
\newcommand{\Fi}{{\bf F}}

\newcommand{\proof}{\mbox{\bf Proof} \hspace{5pt}} 
\newcommand{\remark}{\mbox{\bf Remark} \hspace{5pt}}
\newcommand{\ruimte}{\vskip10.0pt plus 4.0pt minus 6.0pt}

\newcommand{\simh}{{\stackrel{{\rm cap}}{\sim}}}
\newcommand{\ad}{{\mathop{\rm ad}}}
\newcommand{\Ad}{{\mathop{\rm Ad}}}
\newcommand{\Aut}{\mathop{\rm Aut}}
\newcommand{\arccot}{\mathop{\rm arccot}}
\newcommand{\capp}{{\mathop{\rm cap}}}
\newcommand{\rcapp}{{\mathop{\rm rcap}}}
\newcommand{\diam}{\mathop{\rm diam}}
\newcommand{\divv}{\mathop{\rm div}}
\newcommand{\codim}{\mathop{\rm codim}}
\newcommand{\RRe}{\mathop{\rm Re}}
\newcommand{\IIm}{\mathop{\rm Im}}
\newcommand{\Tr}{{\mathop{\rm Tr}}}
\newcommand{\Vol}{{\mathop{\rm Vol}}}
\newcommand{\card}{{\mathop{\rm card}}}
\newcommand{\supp}{\mathop{\rm supp}}
\newcommand{\sgn}{\mathop{\rm sgn}}
\newcommand{\essinf}{\mathop{\rm ess\,inf}}
\newcommand{\esssup}{\mathop{\rm ess\,sup}}
\newcommand{\Int}{\mathop{\rm Int}}
\newcommand{\inn}{\mathop{\rm int}}
\newcommand{\Leibniz}{\mathop{\rm Leibniz}}
\newcommand{\lcm}{\mathop{\rm lcm}}
\newcommand{\loc}{{\rm loc}}

\newcommand{\mod}{\mathop{\rm mod}}
\newcommand{\spann}{\mathop{\rm span}}
\newcommand{\one}{1\hspace{-4.5pt}1}

\hyphenation{groups}
\hyphenation{unitary}

\newcommand{\tfrac}[2]{{\textstyle \frac{#1}{#2}}}

\newcommand{\ca}{{\cal A}}
\newcommand{\cb}{{\cal B}}
\newcommand{\cc}{{\cal C}}
\newcommand{\cd}{{\cal D}}
\newcommand{\ce}{{\cal E}}
\newcommand{\cf}{{\cal F}}
\newcommand{\cg}{{\cal G}}
\newcommand{\ch}{{\cal H}}
\newcommand{\ci}{{\cal I}}
\newcommand{\ck}{{\cal K}}
\newcommand{\cl}{{\cal L}}
\newcommand{\cm}{{\cal M}}
\newcommand{\co}{{\cal O}}
\newcommand{\cs}{{\cal S}}
\newcommand{\ct}{{\cal T}}
\newcommand{\cx}{{\cal X}}
\newcommand{\cy}{{\cal Y}}
\newcommand{\cz}{{\cal Z}}

\newcommand{\wtozp}{W^{1,2}\raisebox{10pt}[0pt][0pt]{\makebox[0pt]{\hspace{-34pt}$\scriptstyle\circ$}}}
\newlength{\hightcharacter}
\newlength{\widthcharacter}
\newcommand{\covsup}[1]{\settowidth{\widthcharacter}{$#1$}\addtolength{\widthcharacter}{-0.15em}\settoheight{\hightcharacter}{$#1$}\addtolength{\hightcharacter}{0.1ex}#1\raisebox{\hightcharacter}[0pt][0pt]{\makebox[0pt]{\hspace{-\widthcharacter}$\scriptstyle\circ$}}}
\newcommand{\cov}[1]{\settowidth{\widthcharacter}{$#1$}\addtolength{\widthcharacter}{-0.15em}\settoheight{\hightcharacter}{$#1$}\addtolength{\hightcharacter}{0.1ex}#1\raisebox{\hightcharacter}{\makebox[0pt]{\hspace{-\widthcharacter}$\scriptstyle\circ$}}}
\newcommand{\scov}[1]{\settowidth{\widthcharacter}{$#1$}\addtolength{\widthcharacter}{-0.15em}\settoheight{\hightcharacter}{$#1$}\addtolength{\hightcharacter}{0.1ex}#1\raisebox{0.7\hightcharacter}{\makebox[0pt]{\hspace{-\widthcharacter}$\scriptstyle\circ$}}}

\newpage

 \thispagestyle{empty}
 
 \begin{center}
 \vspace*{-1.0cm}
\vspace*{1.5cm}

{\Large{\bf On extensions of   }}\\[3mm] 
{\Large{\bf local Dirichlet forms }}  \\[4mm]
\large Derek W. Robinson$^\dag$ \\[1mm]

\normalsize{January  2016}
\end{center}

\vspace{+5mm}

\begin{list}{}{\leftmargin=1.7cm \rightmargin=1.7cm \listparindent=15mm 
   \parsep=0pt}
   \item
{\bf Abstract} $\;$ 
Let $\ce$  be a Dirichlet form on $L_2(X\,;\mu)$ where $(X,\mu)$ is locally compact
$\sigma$-compact measure space.
Assume  $\ce$ is inner regular, i.e.\ regular  in restriction to functions of compact support,
and local in the sense that
 $\ce(\varphi,\psi)=0$ for all $\varphi, \psi\in D(\ce)$ with $\varphi\,\psi=0$.
We  construct  two  Dirichlet forms $\ce_m$ and $\ce_M$ such
that $\ce_m\leq \ce\leq \ce_M$.
These forms are potentially the smallest and largest such Dirichlet forms.
In particular $\ce_m\supseteq \ce_M$,  $(\ce_M)_m=\ce_m$ and $(\ce_m)_M=\ce_M$.
We analyze the family of local, inner regular, Dirichlet forms  $\cf$ which extend
$\ce$ and  satisfy $\ce_m\leq \cf\leq \ce_M$. 
We prove that the latter bounds are valid if and only if $\cf_M=\ce_M$, or $\cf_m=\ce_m$, or $D(\ce_M)$ is an order ideal of $D(\cf)$.
Alternatively the $\cf$ are characterized by $D(\ce_M)\cap L_\infty(X)$ being an algebraic ideal of
$D(\cf)\cap L_\infty(X)$.
As an application we show that if $\ce$ and $\cf$ are strongly local
then the Ariyoshi--Hino set-theoretic distance is the same for each of the forms $\ce$, $\ce_M$ and $\cf$.
If in addition $\ce_m$ is strongly local then it also defines the same distance.
Finally we characterize the uniqueness condition $\ce_M=\ce_m$ by capacity estimates.
\end{list}

\vfill

\noindent AMS Subject Classification: 31C25, 47D07.

\vspace{0.5cm}

\noindent
\begin{tabular}{@{}cl@{\hspace{10mm}}cl}
$ {}^\dag\hspace{-5mm}$&   Mathematical Sciences Institute (CMA)    &  {} &{}\\
  &Australian National University& & {}\\
&Canberra, ACT 0200 && {} \\
  & Australia && {} \\
  &derek.robinson@anu.edu.au
 & &{}\\
\end{tabular}

\newpage

\setcounter{page}{1}

\newpage

\newpage

\section{Introduction}\label{S1}

The standard theory of Dirichlet forms \cite{Fuk}  \cite{BH}  \cite{FOT} encompasses the quadratic forms associated 
with real, positive, symmetric, second-order, elliptic operators on domains in~$\Ri^d$.
The theory also extends to the  non-symmetric situation both real \cite{MR} and complex \cite{Ouh5}.
The Dirichlet form theory, supplemented with  an intrinsic distance, has led to interesting generalizations of many of 
the structural  properties of the elliptic operators and in particular to those related to the underlying Gaussian structure 
(see, for example, \cite{Mosco} \cite{BM} \cite{Stu3, Stu4, Stu5, Stu2}).
These results are based on the analysis of a fixed local, regular,  Dirichlet form.
In contrast our aim is to analyze the family of possible  Dirichlet form extensions of a given local Dirichlet form.
This is a partial analogue of  the classification of the submarkovian extensions of an elliptic operator, 
i.e.\  the extensions which generate submarkovian semigroups.
The extension problem was considered in  \cite{Fuk}, Section~2.2.3, for the Laplacian on a domain in $\Ri^d$ 
and in \cite{FOT} for a  second-order elliptic operator $H$ with smooth coefficients.
It was established that  $H$  has  both  a minimal   and a maximal submarkovian extension (see \cite{FOT}, Theorem~3.3.1).
The minimal extension corresponds to Dirichlet boundary conditions and the maximal extension to Neumann conditions.
Each  submarkovian extension  is  determined by  a Dirichlet form
extension of the  quadratic form $h$ corresponding to $H$.  
These Dirichlet forms are intermediate to the forms of the maximal and minimal extensions with respect to  the 
usual ordering of quadratic forms.
The form ordering is closely related to the ordering of the corresponding submarkovian semigroups.
The semigroup associated with the intermediate form dominates the semigroup corresponding to Dirichlet boundary conditions
and is dominated by the Neumann semigroup (see \cite{Ouh3} and \cite{Ouh5}, Chapters~2 and 4).
The ordering is also related to the lattice ordering of the domains of the Dirichlet forms.
In fact under a mild assumption of regularity the local Dirichlet forms corresponding to semigroups sandwiched between the 
Neumann and Dirichlet semigroups can be completely characterized by generalized Robin boundary conditions \cite{AW, AW2}.
In the sequel we obtain analogues and extensions of these results for local Dirichlet forms.

Let $\ce$ be a   local Dirichlet form  which is inner regular, 
i.e.\ it satisfies the usual  regularity property in restriction to the functions with compact support.
We construct for each  such $\ce$  a maximal 
Dirichlet form $\ce_M$ and a minimal form $\ce_m$.
The maximal form is a restriction of $\ce$ and the minimal  form  an  extension of $\ce_M$.
But $\ce_m$ is  not necessarily an extension of $\ce$.
The submarkovian semigroup $S^m$ corresponding to $\ce_m$ dominates the semigroup $S^M$ corresponding to $\ce_M$.
The construction of both  forms is by `interior approximation'  and corresponds abstractly to the 
imposition of either  Dirichlet or Neuman boundary conditions even in the absence of
an  explicit boundary.
For a precise formulation we refer to Theorem~\ref{tdrf3.10} in Section~\ref{S3}.
The proof of the theorem and the subsequent analysis rely on a number of structural results which 
are summarized  in Section~\ref{S2}.
Then in Section~\ref{S4} we characterize the  local, inner regular, Dirichlet form extensions $\cf$ of $\ce$ 
which satisfy $\ce_m\leq \cf\leq\ce_M$ as the forms with $\cf_M=\ce_M$ or $\cf_m=\ce_m$.
Again $\cf$ is an extension of $\ce_M$ but $\ce_m$ is not necessarily an extension of $\cf$.
The forms $\cf$ can also be characterized by a number of order properties.
In fact $\ce_m\leq \cf\leq \ce_M$ if and only if the domain $D(\ce_M)$ of $\ce_M$ is an order ideal of the domain $D(\cf)$ of $\cf$  (see \cite{Ouh5} Section~2.3)
or if and only if  $B(\ce_M)$, the bounded functions in $D(\ce_M)$, is an algebraic ideal of $B(\cf)$ (see Proposition~\ref{tdrf2}).
The latter algebraic property is fundamental for much of the analysis.
Finally in Section~\ref{S5} we discuss the set-theoretic distance associated with a strongly local Dirichlet form and show that it has
the same values for all the forms $\cf$ satisfying the foregoing conditions.
We also discuss the uniqueness criterion $\ce_m=\ce_M$ and its characterization in terms of capacity conditions.

\section{Preliminaries}\label{S2}

We begin by summarizing some standard and some not so standard elements  of the theory of quadratic  forms which are necessary in the 
subsequent discussion of extensions of Dirichlet forms.
We refer to \cite{Kat1, Kat1a} for a detailed description of the general theory and to \cite{BH} \cite{FOT} for the theory of Dirichlet forms.
Throughout we assume that  $X$ is  a locally compact $\sigma$-compact
metric space and $\mu$ a positive Radon measure with $\supp \mu = X$.
The corresponding real $L_p$-spaces are denoted  by $L_p(X)$.

\smallskip

Let $\ce$  be a positive, densely-defined, quadratic form on $L_2(X)$ with domain $D(\ce)$ and corresponding 
graph norm $\varphi\in D(\ce)\mapsto \|\varphi\|_{D(\ce)}=(\ce(\varphi)+\|\varphi\|_2^2)^{1/2}$.
If $\ce$ is closed with respect to the graph norm then it is the form of a positive self-adjoint operator $H$,  i.e.\
 $D(\ce)=D(H^{1/2})$ and $\ce(\varphi)=\|H^{1/2}\varphi\|_2^2$,   which  generates a self-adjoint contraction semigroup $S$.
If $\ce$ is closable with respect to the graph norm its closure is denoted by $\overline\ce$.
A subspace $D$ of the domain $D(\ce)$ of the closable form $\ce$ is defined to be a core of $\ce$ if $\overline\ce=\overline{\ce|_D}$.
Note that this definition coincides with that in \cite{Kat1} but differs from that of \cite{FOT}.

\smallskip

The form $\ce$ is defined to be Markovian if for each $\varphi\in D(\ce)$ one has $ \varphi\wedge 1\in D(\ce)$ and 
\begin{equation}
\ce( \varphi\wedge 1)\leq \ce(\varphi)
\label{edrf2.1}
\end{equation}
and $\ce$ is defined to be Dirichlet if it is both closed and Markovian.
The Markovian property extends to a much broader class of mappings.
A map $F\colon \Ri\mapsto \Ri$ is called a normal contraction if $|F(x)-F(y)|\leq |x-y|$ and $F(0)=0$.
If $\ce$ is a Dirichlet form and $F$ a normal contraction then $F\circ(D(\ce))\subseteq D(\ce)$ and
\begin{equation}
\ce(F\circ\varphi)\leq \ce(\varphi)
\label{edrf2.11}
\end{equation}
for all $\varphi\in D(\ce)$. 

\smallskip

If $\ce$ is a Dirichlet form the corresponding semigroup $S$ is submarkovian, i.e.\ if $0\leq \varphi\leq \one$ then $0\leq S_t\varphi\leq \one$
for all $t>0$.
Conversely if $S$ is a self-adjoint submarkovian semigroup then the quadratic form $\ce$ corresponding to its generator $H$  is Dirichlet. 

If $\ce$ is Dirichlet the subspace $B(\ce)=D(\ce)\cap L_\infty(X)$ of bounded functions in the domain $D(\ce)$
is an algebra and also a core of $\ce$.
(See  \cite{BH} Sections~1.1--1.3).
We let $B_c(\ce)$ denote the subalgebra of $B(\ce)$ spanned by the functions with compact support  and set $C_c(\ce)=D(\ce)\cap C_c(X)=B_c(\ce)\cap C(X) $.

\smallskip

In the sequel we consider two regularity properties for the Dirichlet form $\ce$.
These properties are defined with the aid of the following conditions.
 \begin{equation}
 \left. \begin{array}{ll}
1.&C_c(\ce) \mbox{ is dense in } C_0(X) \mbox{ with respect to the supremum norm,}\hspace{1cm}\\[8pt]
2.&C_c(\ce) \mbox{ is dense in } B_c(\ce) \mbox{ with respect to the graph norm},
  \end{array}\right\}\label{edrf2.00}
\end{equation}
where $C_0(X)$ denotes the space of continuous functions over $X$ which vanish at infinity.
(This notation again differs from that of \cite{FOT}.)
We define $\ce$ to be semi-regular if Condition~1 of (\ref{edrf2.00}) is satisfied and inner regular if both Conditions~1 and 2 are satisfied.
The standard definition of regularity (see \cite{FOT} Section~1.1) replaces Condition~2 with the stronger assumption that $C_c(\ce)$
is dense in $B(\ce)$ or, equivalently, dense in $D(\ce)$ but the latter assumption is not widely applicable 
in the analysis of extensions of Dirichlet forms.
This  is    illustrated by  the Dirichlet forms associated with  the Laplacian on a domain $\Omega$  in Euclidean space $\Ri^d$.
The form corresponding to Dirichlet boundary conditions is given by $\ce_D(\varphi)=\|\nabla\varphi\|_2^2$ with $D(\ce_D)=W^{1,2}_{\,0}(\Omega)$ and 
it is regular in the sense of \cite{FOT}.
But the form corresponding to Neumann boundary conditions,  $\ce_N(\varphi)=\|\nabla\varphi\|_2^2$ with $D(\ce_N)=W^{1,2}(\Omega)$, is inner regular
although  it is not regular
unless, exceptionally, $W^{1,2}_{\,0}(\Omega)=W^{1,2}(\Omega)$.
Regularity fails, for the same reason, for  extensions with  mixtures of Dirichlet and Neumann conditions or  with Robin boundary conditions
although inner regularity persists.

Semi-regularity has the following local implications.
\begin{lemma}\label{l2}
Assume  $\ce$ is semi-regular.  
Let $Y$ be a bounded open subset of $X$.
\begin{tabel}
\item\label{l2-1}
$D(\ce)\cap C_c(Y)$ is dense in $C_0(Y)$ with respect to the supremum norm and 
in $L_2(Y)$ with respect to the $L_2$-norm.
\item\label{l2-2}
There is an $\eta\in C_c(\ce)$  with $0\leq \eta\leq 1$ and $\eta=1$ on $Y$.
\end{tabel}
\end{lemma}
\proof\  \ref{l2-1}. 
It suffices to establish the density in $C_0(Y)$ since this implies
 density in $L_2(Y)$.
But each $\varphi\in C_0(Y)$ has a unique decomposition  into a positive and a negative component.
Therefore  it
suffices to prove that each positive $\varphi\in C_0(Y)$ can be approximated uniformly by a sequence of $\varphi_n\in D(\ce)\cap \,C_c(Y)$.
It follows, however,  from the semi-regularity that for each $n\in\Ni$ there is a $\psi_n\in C_c(\ce)$ such that $\|\varphi-\psi_n\|_\infty\leq 1/(2n)$.
Replacing $\psi_n$ by $\psi_n\wedge0$ does not affect this estimate so one can assume that $\psi_n$ is positive.
Then one must have $0\leq \psi_n(y)\leq 1/(2n)$ for all $y\in Y^{\rm c}$.
Hence  $\varphi_n=\psi_n-(\psi_n\wedge(1/2n)) \in D_Y(\ce)$ and $\|\varphi-\varphi_n\|_\infty<1/n$.
\smallskip

\noindent \ref{l2-2}. There is a $\varphi\in C_0(X)$ such that $\varphi\geq 2$ on $Y$.
But it follows from semi-regularity of $\ce$ that there is a $\psi\in C_c(\ce)$ such that
$\|\varphi-\psi\|_\infty\leq1$.
Therefore $\psi\geq1$ on $Y$.
Then $\eta=0\vee(\psi\wedge1)\in C_c(\ce)$ and $\eta=1$ on $Y$.\hfill$\Box$

\bigskip

One can associate with each  Dirichlet form $\ce$ and each  $\xi\in B(\ce)$ the truncated form $\ce_\xi$ 
by setting  $D(\ce_\xi)=B(\ce)$ and 
\begin{equation}
\ce_\xi(\varphi,\psi)= 2^{-1}( \ce(\xi\,\varphi, \psi)+ \ce(\varphi, \xi\,\psi) - \ce(\xi,\varphi\,\psi)
\label{edrf2.20}
\end{equation}
for all $\varphi,\psi\in B(\ce)$. 
Then
\begin{equation}
\ce_\xi(\varphi)=\ce_\xi(\varphi,\varphi)=\ce(\varphi,\xi\,\varphi)-2^{-1}\ce(\xi,\varphi^2)
\label{edrf2.2}
\end{equation}
for all $\varphi\in B(\ce)$. 
If $\xi\in B(\ce)_+$ then $\ce_\xi$ is a Markovian form  which satisfies  the bounds $0\leq \ce_\xi(\varphi) \leq \|\xi\|_\infty \, \ce(\varphi)$ for all $\varphi\in B(\ce)$
(see \cite{BH}, Proposition~I.4.1.1).
Therefore the quadratic form $\ce_\xi$ extends by continuity to a Markovian form on $D(\ce)$ although the identity (\ref{edrf2.2}) is not necessarily valid for the extension.
Moreover, if $\ce$ is semi-regular the definition of the truncated forms can be extended to  all $\xi\in C_0(X)$.
Then $\xi\in C_0(X)\mapsto \ce_\xi(\varphi)$ is a positive linear functional for each $\varphi\in D(\ce)$.
Therefore there is a Radon measure $\mu_\varphi$ such that $\mu_\varphi(\xi)=\ce_\xi(\varphi)$ for all $\xi\in C_0(X)$.

\smallskip

Although the truncated forms $\ce_\xi$ are Markovian they are not necessarily closed nor closable.
Therefore we consider  their relaxations.
The  relaxation  of a positive quadratic form is defined as 
the largest positive closed form $\ce_0$ which is dominated by $\ce$, i.e.\ the largest positive  closed form $\ce_0$ such that $D(\ce)\subseteq D(\ce_0)$ and $\ce_0(\varphi)\leq \ce(\varphi)$ for all $\varphi\in D(\ce)$.
The  relaxation is also  called the  lower semi-continuous regularization (see \cite{ET}, page~10) or the  relaxed form (see \cite{DalM}, page~28).
Alternatively, it can be characterized as the closure of the regular part of the form (see \cite{bSim5}).
Note that  $\ce_0(\varphi)=\ce(\varphi)$ for all $\varphi\in D(\ce)$ if and only if $\ce$ is closable and  then $\ce_0=\overline \ce$.

The domain  $D(\ce_0)$  of the relaxation is the subspace of $L_2(X)$
spanned by the $\varphi\in L_2(X)$ for which there is a sequence $\{\psi_n\}_{n\geq1}$  with $\psi_n\in D(\ce)$ such that $\lim_{n\to\infty}\|\psi_n-\varphi\|_2=0$
and $\liminf_{n\to\infty}\ce(\psi_n)<\infty$.
Moreover,
\begin{equation}
\ce_0(\varphi)=\liminf_{\psi\to\varphi} \ce(\psi)
\label{edrf2.3}
 \end{equation}
where $\liminf_{\psi\to\varphi} \ce(\psi)$ indicates the infimum over  the  $\liminf_{n\to\infty}\ce(\psi_n)$ for all choices of $\psi_n\in D(\ce)$ such that $\lim_{n\to\infty}\|\psi_n-\varphi\|_2=0$ (see  \cite{bSim4}, Theorem~3).
If $\ce$ is closable then $\ce_0=\overline\ce$ and this characterization follows from
 \cite{Kat1}, Theorem~VI.1.16.
 In analogy with the closable case we define a subspace $D$ of $D(\ce)$ to be a core of $\ce$ if $\ce_0=(\ce|_D)_0$.

The  relaxation procedure has  many properties similar to the closure.
In particular if the two forms $\ce$ and $\cf$  satisfy $0\leq\ce\leq\cf$ then  their relaxations satisfy $0\leq \ce_0\leq \cf_0$.
The following observation is used regularly.

\begin{lemma}\label{ldrf2.2}
The relaxation of a Markov form is a Dirichlet form.
\end{lemma}
\proof\
Let $\ce_0$ denote the relaxation of the Markov form $\ce$.
If    $\varphi \in D(\ce_0)$  then  there are $\psi_n \in D(\ce)$ such that 
$\lim_{n \to \infty} \| \psi_n - \varphi\|_2=0$ and 
$ \lim_{n \to \infty}  \ce(\psi_n)<\infty$.
But  since $x\mapsto x\wedge 1$ is a normal contraction (see \cite{BH}, Section~I.2.3) it follows that  $\lim_{n \to \infty}  \|\psi_n \wedge \one-\varphi \wedge \one\|_2=0$. 
Moreover,  $\liminf_{n \to \infty}  \ce( \psi_n \wedge \one) \leq \liminf_{n \to \infty} \ce(\psi_n)<\infty$ since $\ce$ is Markovian.
This, however,  implies that  $ \varphi \wedge \one \in D(\ce_0)$ and $\ce_0(\varphi \wedge \one) \leq \ce_0(\varphi)$.
Thus $\ce_0$ is closed and Markovian, i.e.\ Dirichlet.
\hfill$\Box$

\bigskip

Much of the subsequent analysis is for local Dirichlet forms but there are various notions of locality. 
In particular the definitions of \cite{BH} and \cite{FOT} are different.
We  adopt  definitions  intermediate between these two texts.
Specifically we define the  Dirichlet form $\ce$  to be local if  $\ce(\varphi,\psi)=0$ for all $\varphi,\psi\in D(\ce)$   such that  $\varphi\,\psi=0$.
Further we define $\ce$ to be strongly local if  $\ce(\varphi,\psi)=0$ for all $\varphi,\psi\in D(\ce)$   such that  $(\varphi+a\one)\,\psi=0$ for some  $a\in\Ri$.
These definitions are  similar to those of  \cite{FOT}, Section~1.1, where they are only required for functions of compact support.
The second condition coincides exactly with the definition of locality in  \cite{BH}, Section~I.1.5.
In fact the locality conditions are only necessary for all pairs $\varphi,\psi$ in a suitable core.

\begin{prop}\label{pdrf2.3}
Let $D$ be a core of the Dirichlet form $\ce$.
Assume $\varphi\in D$ implies $|\varphi|\in D$.
Then the  following conditions are equivalent:
\begin{tabelpairs}
\item\label{pdrf2.3-1}
 $\ce(\varphi,\psi)=0$ for all $\varphi,\psi\in D(\ce)$   $($for all $\varphi,\psi\in D)$ such that  $\varphi\,\psi=0$,
 \item\label{pdrf2.3-3}
$ \ce(|\varphi|)=\ce(\varphi)$ for all $\varphi\in D(\ce)$ $($for all $\varphi\in D)$.
\end{tabelpairs}
\end{prop}
\proof\ Clearly \ref{pdrf2.3-1}$\Rightarrow$\ref{pdrf2.3-1}$^\prime$.

\smallskip

\noindent \ref{pdrf2.3-1}$^\prime\!\!\Rightarrow$\ref{pdrf2.3-3}$^\prime.\;$
If $\varphi\in D$ then $|\varphi|\in D$ and $\varphi_\pm=(|\varphi|\pm\varphi)/2\in D$.
But $\varphi_+\,\varphi_-=0$. 
Therefore Condition~\ref{pdrf2.3-1}$^\prime$ gives
$\ce(|\varphi|)=\ce(\varphi_++\varphi_-)=\ce(\varphi_+)+\ce(\varphi_-)=\ce(\varphi_+-\varphi_-)=\ce(\varphi)$.
\smallskip

\noindent \ref{pdrf2.3-3}$^\prime\!\!\Rightarrow$\ref{pdrf2.3-3}.$\;$
This follows from Ancona's proof, \cite{Anc} Proposition~4, of the continuity of the  contraction $\varphi\mapsto |\varphi|$  in the graph norm of $D(\ce)$.
In fact Theorem~10 of \cite{Anc} establishes that all the normal contractions acting on  $D(\ce)$ are continuous.  

\smallskip

\noindent \ref{pdrf2.3-3}$\Rightarrow$\ref{pdrf2.3-1}.$\;$
Let $\varphi, \psi\in D(\ce)_+$ with $\varphi\,\psi=0$.
Set $\chi=\varphi-\psi$ then $|\chi|=\varphi+\psi$ and $\ce(|\chi|)-\ce(\chi)=2\,\ce(\varphi,\psi)$.
Therefore $\ce(\varphi,\psi)=0$ by Condition~\ref{pdrf2.3-3}.
This establishes the locality condition on positive elements of $D(\ce)$.
Next assume $\varphi, \psi\in D(\ce)$ with $\varphi\,\psi=0$ and with  positive
and negative components  $\varphi_\pm, \psi_\pm$.
Then $\varphi_+\,\varphi_-=0$ and $\psi_+\,\psi_-=0$.
Therefore $(\varphi_\pm\,\psi_\pm)^2=\varphi_\pm\,(\varphi\,\psi)\,\psi_\pm=0$.
Thus all  products of the $\varphi_\pm$ with the $\psi_\pm$ are zero.
Then $\ce(\varphi,\psi)=0$  by linearity from the previous conclusion for positive elements.
\hfill$\Box$

\bigskip

There is a similar characterization of strong locality but the mapping $x\mapsto |x|$ is replaced by the normal contraction $x\mapsto F(x)=|x+1|-1$, or by $x\mapsto F(x)=|x+\lambda|-\lambda$ with $\lambda>0$.
The details of the strongly local case are rather different to the previous one and depend on yet another characterization of strong locality
(see \cite{BH} Section~I.5 and in particular Proposition~I.5.1.3 and Remark~I.5.1.5).

\smallskip

It follows from Proposition~\ref{pdrf2.3} that if 
$\ce$ is a  local, closable, Markov form whose domain is invariant under the map $\varphi\mapsto |\varphi|$ then its closure is a local Dirichlet form.
It does not, however, automatically follow that  a similar conclusion is valid for the relaxation of a  local Markovian form which is not closable
(see \cite{Mosco}, Example~6.1.1).
There are, however, positive results of this type \cite{FOT}, Theorem~3.1.2,  and \cite{ERSZ2}, Proposition~2.2.

\smallskip

The subsequent  discussion of upper and lower bounds for local Dirichlet forms depends in part on two monotonic approximation arguments.

\begin{prop}\label{pdrf2.4}
Let $\ce_\alpha$ be a monotonically increasing net of positive, closed,  quadratic forms.
Define the limit form $\ce$
by $D(\ce)=\{\varphi: \varphi\in\bigcap_\alpha D(\ce_\alpha)\,, \sup_\alpha \ce_\alpha(\varphi)<\infty\}$ and 
\[
\ce(\varphi)=\lim_\alpha\ce_\alpha(\varphi)=\sup_\alpha\ce_\alpha(\varphi)
\;.
\]
Then $\ce$ is positive and closed.
If $\ce$, and consequently all the $\ce_\alpha$, are densely-defined then
the positive self-adjoint operators $H_\alpha$ corresponding to the $\ce_\alpha$ converge in the strong resolvent sense to the
positive self-adjoint operator $H$ corresponding to $\ce$.
\end{prop}

This result was substantially established by \cite{Kat1}, Theorem~VIII.3.13, but this early
 version assumed  that the limit form was $\ce$ closed.
Subsequently it was observed  in  \cite{Robl}, Section~I.2.9, that $\ce$ is automatically closed (see  also \cite{BR2}, Lemma~5.2.13). 
Somewhat later Kato reached  the same conclusion by quite different lower semicontinuity arguments (see \cite{Kat1a} Theorem~VIII.3.13a) and Simon \cite{bSim5} gave a third distinct proof.
It   also follows that if the $\ce_\alpha$ are closable then the limit form is closable since the monotonic limit of the closures $\overline{\ce}_{\!\alpha}$ is a closed extension
of the monotonic limit of the $\ce_\alpha$.

\begin{prop}\label{pdrf2.5}
Let $\ce_\alpha$ be a monotonically decreasing  net of positive, closed, densely-defined, quadratic forms.
Define $\ce$
by $D(\ce)=\{\varphi: \varphi\in\bigcup_\alpha D(\ce_\alpha)\}$ and 
\[
\ce(\varphi)=\lim_\alpha\ce_\alpha(\varphi)=\inf_\alpha\ce_\alpha(\varphi)
\;.
\]
Then $\ce$ is positive and densely-defined.
The positive self-adjoint operators $H_\alpha$ corresponding to the $\ce_\alpha$ converge in the strong resolvent sense to the
positive self-adjoint operator $H_0$ corresponding to the  relaxation $\ce_0$ of $\ce$.
\end{prop}

A version of this result was given by Kato, \cite{Kat1} Theorem~VIII.3.11, again with the assumption that the limit form is closable. 
The complete statement was derived by   Simon \cite{bSim5} with a slightly different terminology.

Propositions~\ref{pdrf2.4} and \ref{pdrf2.5} also have extensions to monotone families of non-densely defined forms,
e.g.\ forms $\ce_\alpha$ defined on subspaces $L_2(X_\alpha)$ of $L_2(X)$ (see \cite{bSim5}, Section~4).
This observation will be used in the subsequent discussion of order properties of submarkovian semigroups
(see Propositions~\ref{tdrf2}, \ref{pdrf3.11}, \ref{pdrf3.13} and \ref{pdrf3.14}).
Note that in both propositions, and their extensions, the strong resolvent convergence of the semigroup generators is equivalent to the  strong convergence
of the semigroups
 $S^\alpha$ generated by the $H_\alpha$  to the semigroup $S$ generated by the limit operator $H$ or $H_0$ (see \cite{Kat1}, Theorem~IX.2.16).
 
 \smallskip
 
 A measurable  subset $A$ of $X$ is defined to have   finite capacity with respect to the Dirichlet form $\ce$, or finite $\ce$-capacity,
if  there  is a $\varphi\in B(\ce)$ with  $\varphi= 1$ on $A$.
In particular if $\ce$ is semi-regular then  each bounded open subset
has finite capacity.
The finite capacity subspace of $D(\ce)$  is defined by
\[
D_{\capp}(\ce)=\{\varphi\in B(\ce): \supp\varphi \mbox{ has finite capacity }\}
\;.
\]
It is  invariant under normal contractions.
Set $B_{\capp}(\ce)=D_{\capp}(\ce)\cap L_\infty(X)$.

\begin{prop}\label{pcap1} 
The subspace $B_{\capp}(\ce)$ is a core of $\ce$.
\end{prop}
\proof\
Since $B(\ce)$ is a core it suffices to prove that each $\varphi\in B(\ce)_+$  can be approximated in the $D(\ce)$-graph norm by
a sequence $\psi_n\in B_{\capp}(\ce)$.
Define $\varphi_\tau$ by $\varphi_\tau(x)=\tau^{-1}(\varphi(x)\wedge \tau)$ with $\tau>0$.
Then $\varphi_\tau\in B(\ce)$, $0\leq \varphi_\tau\leq 1$ 
and $\ce(\varphi_\tau)\leq \tau^{-2}\ce(\varphi)$ by the Dirichlet property.
Therefore the set $A_\tau=\{x\in X: \varphi_\tau(x)=1\}=\{x\in X: \varphi(x)\geq \tau\}$ has finite capacity.
Now set $\psi_\sigma=\varphi -\varphi\wedge \sigma^{-1}$ with $\sigma>0$.
Then  $\psi_\sigma\in B(\ce)$ and $\supp\psi_\sigma=A_{\sigma^{-1}}$.
Therefore $\psi_\sigma\in B_{\capp}(\ce)$.
But $\|\varphi-\psi_\sigma\|_{D(\ce)}=\|\varphi\wedge \sigma^{-1}\|_{D(\ce)}\to0$ as $\sigma\to\infty$ by the Dirichlet form structure
(see, for example, \cite{FOT} Theorem~1.4.2(iv)).
This establishes the required approximation.\hfill$\Box$

\bigskip

Finally we consider various  various ordering properties of two Dirichlet forms $\ce$, $\cf$  and the corresponding submarkovian semigroups $S$, $T$.
The ordering of the semigroups, i.e.\ the property $0\leq S_t\varphi\leq T_t\varphi$ for all $\varphi\in L_2(X)_+$ and $t>0$, is related
 to an order ideal property of the form domains (see \cite{Ouh5}, Theorem~2.24).
This observation is independent of the detailed Dirichlet structure and depends only on the positivity
and $L_2$-contractivity of the semigroups.
But for Dirichlet forms the ordering is  also equivalent
to an algebraic ideal property.

First if $\ce$ and $\cf$ are Dirichlet forms then Ouhabaz defines $D(\ce)$ to be an ideal of $D(\cf)$ if
 $\varphi\in D(\ce)$, $\psi\in D(\cf)$  and $|\psi|\leq |\varphi|$  implies that $\psi\sgn\varphi\in D(\ce)$.
Alternatively $D(\ce)$ is defined to be an order ideal of $D(\cf)$ if  $0\leq\psi\leq \varphi$ with $\varphi\in D(\ce)$ and $\psi\in D(\cf)$
implies that $\psi\in D(\ce)$.
These two notions are closely related and if $D(\ce)\subseteq D(\cf)$ it follows that they are equivalent,
i.e.\ $D(\ce)$ is an ideal of $D(\cf)$ if and only if it is an order ideal.
This follows from the proof of Proposition~2.23 in \cite{Ouh5}.

Secondly, since $\ce$ and $\cf$ are Dirichlet forms $B(\ce)$ and $B(\cf)$ are algebras and  $B(\ce)$ is defined to be an algebraic ideal of $B(\cf)$ if $B(\ce)\,B(\cf)\subseteq B(\ce)$.
This property is also related to the ordering of the corresponding semigroups $S$ and $T$.

\begin{prop}\label{tdrf2}
Let $\ce$ and $\cf$ be  Dirichlet forms with  $D(\ce)\subseteq D(\cf)$ and 
$S$,  $T$  the submarkovian semigroups associated with $\ce$ and $\cf$, respectively.

The following conditions are equivalent:
\begin{tabel}
\item\label{tdrf2-1}
$0\leq S_t\varphi\leq T_t\varphi$ for all $\varphi\in L_2(X)_+$ and all $t>0$,
\item\label{tdrf2-2}
$ D(\ce)$ is an order ideal of $ D(\cf)$ and $\ce(\varphi,\psi)\geq \cf(\varphi,\psi)$ for all $\varphi,\psi\in D(\ce)_+$,
\item\label{tdrf2-3}
$B(\ce)$ is an algebraic ideal of $B(\cf)$  and $\ce(\varphi,\psi)\geq \cf(\varphi,\psi)$ for all $\varphi,\psi\in D(\ce)_+$.
\end{tabel}
\end{prop}
\proof\
\ref{tdrf2-1}$\Leftrightarrow$\ref{tdrf2-2}.$\;$ 
This equivalence follows from Corollary~2.22 and Proposition~2.23 of \cite{Ouh5}.
In particular the off-diagonal order property of the forms follows from the order property of the semigroups since
\[
\ce(\varphi,\psi)=\lim_{t\to0}t^{-1}(\varphi, (I-S_t)\psi)\geq\lim_{t\to0}t^{-1}(\varphi, (I-T_t)\psi)= \cf(\varphi,\psi)
\]
for all $\varphi,\psi\in D(\ce)_+$.

\smallskip

\noindent\ref{tdrf2-2}$\Rightarrow$\ref{tdrf2-3}.$\;$ 
If $\varphi\in B(\ce)_+$ and $\psi\in B( \cf)_+$ then 
 $ 0\leq \varphi\,\psi\leq \|\psi\|_\infty\,\varphi\in B(\ce)$.
Therefore $\varphi\,\psi\in  B(\ce)$ by  the order ideal property.
 This establishes the algebraic ideal property for positive functions and it then 
 follows for general  $\varphi\in B(\ce)$ and $\psi\in B( \cf)$ 
 by decomposition into positive and negative components.

\smallskip

\noindent\ref{tdrf2-3}$\Rightarrow$\ref{tdrf2-2}.$\;$ 
If $0\leq\psi\leq \varphi$ with  $\varphi\in B(\ce)$ and $\psi\in B(\cf)$  
then one must prove that $\psi\in B(\ce)$.

First choose $\chi_n\in B_{{\rm cap}}(\ce)$ such that $\|\chi_n-\varphi\|_{D(\ce)}\to0$ as $n\to\infty$.
This is possible by Proposition~\ref{pcap1}.
Next the modulus map $x\to|x|$ is a normal contraction and the normal contractions acting on $D(\ce)$ are strongly continuous
with respect to the $D(\ce)$-graph norm (see \cite{Anc}, Proposition~4 and Theorem~10).
Hence  $\||\chi_n|-\varphi\|_{D(\ce)}\to0$ as $n\to\infty$.
Now if $\varphi_n=|\chi_n|\wedge\varphi$ then $\varphi_n\in B_{{\rm cap}}(\ce)$ since $\supp\varphi_n\subseteq \supp\chi_n$, $0\leq \varphi_n\leq \varphi$.
Further
\[
\|\varphi_n-\varphi\|_{D(\ce)}=\|(\varphi-|\chi_n|)_+\|_{D(\ce)}
\]
by the lattice relation $x\wedge y=y-(y-x)_+$.
But $F(x)=x_+$ is a normal contraction.
Therefore $\|(\varphi-|\chi_n|)_+\|_{D(\ce)}\to0$ as $n\to\infty$  again by continuity of the normal contractions.
Hence $\varphi_n\in B_{{\rm cap}}(\ce)$, $0\leq \varphi_n\leq \varphi$ and $\|\varphi_n-\varphi\|_{D(\ce)}\to0$.

Secondly set $\psi_n=\varphi_n\wedge\psi$.
Then $ B_{{\rm cap}}(\ce)\subseteq B_{{\rm cap}}(\cf)$ by the assumption $D(\ce)\subseteq D(\cf)$.
Therefore $\varphi_n\in B_{{\rm cap}}(\ce)\subseteq B_{{\rm cap}}(\cf)$ and $\psi\in B(\cf)_+$  from which it follows that
$\psi_n\in B_{{\rm cap}}(\cf)$.
Moreover, because $\varphi_n\in B_{{\rm cap}}(\ce)$
one can choose $\eta_n\in B(\ce)$ with $0\leq\eta_n\leq1$ and $\eta_n=1$ on $\supp\varphi_n$.
But  $\supp\psi_n\subseteq \supp\varphi_n$.
Therefore
\[
\psi_n=\eta_n\,\psi_n\in B(\ce)\,B(\cf)\subseteq B(\ce)
\]
by the algebraic ideal property.
Then 
\[
\|\psi_n-\psi\|_{D(\cf)}=\|(\psi-\varphi_n)_+\|_{D(\cf)}
\]
by another application of the lattice relation.
Now $\|\varphi_n-\varphi\|_{D(\ce)}=\|\varphi_n-\varphi\|_{D(\cf)}\to0$ as $n\to\infty$.
Therefore $\|(\psi-\varphi_n)_+\|_{D(\cf)}\to  \|(\psi-\varphi)_+\|_{D(\cf)}$
as $n\to\infty$ by continuity of the normal contractions.
But  $(\psi-\varphi)_+=0$ since  $0\leq\psi\leq\varphi$.
Hence $\|(\psi-\varphi_n)_+\|_{D(\cf)}\to0$ as $n\to\infty$.
Therefore $\|\psi_n-\psi\|_{D(\cf)}\to0$.
Finally as the $\psi_n\in B(\ce)$ it follows that they form a Cauchy sequence with respect
to the $D(\ce)$-graph norm.
Consequently $\psi\in B(\ce)$.\hfill$\Box$

\bigskip

The assumption $D(\ce)\subseteq D(\cf)$ in the proposition follows immediately if $\ce\subseteq \cf$.
But in this latter case the
off-diagonal bounds  $\ce(\varphi,\psi)\geq \cf(\varphi,\psi)$ for  $\varphi,\psi\in D(\ce)_+$ are obvious.
Alternatively  $D(\ce)\subseteq D(\cf)$ is a consequence of the order relation $\cf\leq \ce$ but the off-diagonal
bounds are not immediate from this latter relation.
In conclusion we note that the  foregoing proposition remains valid if $\ce$ is not densely-defined  but is a Dirichlet form on a  closed subspace
of $L_2(X)$ (see \cite{Ouh5}, Section~2.6).

\section{Extremal forms}\label{S3}
Let $\ce$  denote a local inner regular Dirichlet form on $L_2(X)$.
In this section we construct two auxiliary Dirichlet forms $\ce_m$ and $\ce_M$ with the property
that $\ce_m\leq \ce\leq \ce_M$.
The construction is such that the forms are extremal in a sense to be established  in the 
following section.
Although the principal results require inner regularity of $\ce$ the construction of the forms $\ce_m$, $\ce_M$ and a number of the intermediate 
statements only require semi-regularity.

First define $\ce_M$ as the closure of the restriction of $\ce$ to $C_c(\ce)$.
Then, by semi-regularity of $\ce$, the form $\ce_M$   is densely-defined.
It  is automatically a regular Dirichlet form which bounds $\ce$ from above, i.e.\ $\ce\leq \ce_M$.
Note that  $C_c(\ce_M)=C_c(\ce)$ and  $\ce=\ce_M$ if and only if $\ce$ is regular.
It also follows from inner regularity that $\ce_M$ is the closure of $\ce$ restricted to $B_c(\ce)$.
In fact the converse is valid.
The closure of the restriction of $\ce$ to $C_c(\ce)$ is equal to the closure of the restriction to $B_c(\ce)$
if and only if $\ce$ is inner regular.

Secondly, to construct the lower bound $\ce_m$ we introduce  the family
$\cy$ of bounded open subsets $Y$ of $X$.
The family  $\cy$ becomes a directed set when the subsets~$Y$ are ordered by inclusion.
For each $Y\in\cy$  define  the convex subset $\cc_Y(\ce)$ of $C_c(\ce)$ by
\[
\cc_Y(\ce)
=\{\xi\in C_c(\ce)\,:\,0\leq\xi\leq \one_Y\}
\;.
\]
Then  $\cc_Y(\ce)$ is also a directed set with respect to the natural order, e.g.\ if $\xi_1,\xi_2\in \cc_Y(\ce)$
then $\xi_{12}=\xi_1\xi_2-\xi_1-\xi_2\in \cc_Y(\ce)$ is a common upper bound.
But   the truncated forms 
$\{\ce_\xi\}_{\xi\in\cc_Y(\ce)}$ 
form a monotonically increasing net of Markovian forms which are uniformly bounded by $\ce$.
 Therefore  one can define a  Markovian  form $\ce_{m,Y}$ by $D(\ce_{m,Y})=B(\ce)$ and 
\begin{equation}
\ce_{m,Y}(\varphi)
=\lim_{\xi\in\cc_Y(\ce)}\ce_{ \xi}(\varphi)
=\sup_{\xi\in\cc_Y(\ce)}  \ce_{ \xi}(\varphi)
\label{edrf3.000}
\end{equation}
for all $\varphi\in B(\ce)$.
Since  $\ce_{\xi}\leq \ce$ for all $\xi\in\cc_Y(\ce)$ one has
$\ce_{m,Y}\leq \ce$ and $\ce_{m,Y}$ extends to $D(\ce)$ by continuity.
Although the forms $\ce_{m,Y}$ are  Markovian they are not necessarily
closed nor even closable.
Nevertheless their relaxations $\ce_{m,Y;\,0}$ are Dirichlet forms on $L_2(X)$ by Lemma~\ref{ldrf2.2}.
Next $\ce_\xi\leq \ce_\eta$ for $0\leq \xi\leq \eta$ and $\cc_{Y_1}(\ce)\subseteq \cc_{Y_2}(\ce)$ for all $Y_1,Y_2\in\cy$  with $Y_1\subseteq Y_2$.
Therefore $\ce_{m,Y_1}\leq \ce_{m,Y_2}\leq \ce$.
Consequently  $\ce_{m,Y_1;\,0}\leq \ce_{m,Y_2;\,0}\leq \ce$.
Thus the relaxations $\{ \ce_{m,Y;\,0}\}_{Y\in \cy}$ are a uniformly bounded monotonically increasing net of Dirichlet forms.
One can then define a Dirichlet form $\ce_m$ on $L_2(X)$ by 
\begin{equation}
D(\ce_m)=\{\varphi\in\bigcap_{Y\in \cy}D(\ce_{m,Y;\,0}):\sup_{Y\in\cy}\ce_{m,Y;\,0}(\varphi)<\infty\}
\label{edfr3.00}
\end{equation}
and  
\begin{equation}
\ce_m(\varphi)=\lim_{Y\in \cy}\ce_{m,Y;\,0}(\varphi)=\sup_{Y\in\cy}\ce_{m,Y;\,0}(\varphi)
\label{edfr3.0}
\end{equation}
for all $\varphi\in D(\ce_m)$.
It follows immediately that $D(\ce)\subseteq D(\ce_m)$ and $\ce_m\leq \ce$.
Note, however,  that $\ce_m$ is not necessarily an extension of $\ce$.

Locality of  $\ce$ ensures that each $\ce_{m,Y}$  is also local.
Moreover,  $\ce_{m,Y}$  is localized on $Y$. 
In particular if $\varphi\in B(\ce)$ and $\supp\varphi\in Y^{\rm c}$ then $\xi\,\varphi=0$ for all
$\xi\in \cc_Y(\ce)$ and it follows that $\ce_{m,Y}(\varphi)=0$.
Moreover, if $\varphi,\psi\in B(\ce)$  and $\varphi|_Y=\psi|_Y$ then it follows from locality that  $\ce_{m,Y}(\varphi)=\ce_{m,Y}(\psi)$.

As an illustration let $\Omega$ be a bounded open subset $\Ri^d$  and consider the Dirichlet form $\ce$ defined as the closure of $\|\nabla\varphi\|_2^2$
on the domain  $\{\varphi|_\Omega:\varphi\in C_c^\infty(\Ri^d\backslash F)\}$ where  $F$ is a closed subset of the boundary $\partial\Omega$.
The corresponding self-adjoint operator is the version of the Laplacian on $L_2(\Omega)$ with Dirichlet boundary conditions on $F$
and Neumann on $\partial\Omega\backslash F$.
It follows immediately from the foregoing construction that $\ce_M(\varphi)=\|\nabla\varphi\|_2^2$  with domain $W^{1,2}_{\,0}(\Omega)$ and $\ce_m\,(=\ce_{m,\Omega})$ is the extension of $\ce_M$
to the domain $W^{1,2}(\Omega)$.
The maximal form corresponds to Dirichlet  conditions on $\partial\Omega$ and the minimal form to Neumann conditions.

\smallskip

\begin{thm}\label{tdrf3.10}
Assume $\ce$ is a  local, inner regular, Dirichlet form.
Then the following following properties are valid:
\begin{tabel}
\item\label{tdrf3.10-1}
$(\ce_M)_m= \ce_m$,
\item\label{tdrf3.10-2}
$(\ce_m)_M=\ce_M$,
\item\label{tdrf3.10-3}
$\ce_m\supseteq \ce_M$,
\item\label{tdrf3.10-4}
  $B_c(\ce)=B_c(\ce_m)=B_c( \ce_M)$,
\item\label{tdrf3.10-41}
$C_c(\ce)=C_c(\ce_m)=C_c( \ce_M)$,
\item\label{tdrf3.10-5}
$\ce_M$ is regular and  $\ce_m$  is inner regular.
   \end{tabel}
 \end{thm}
\proof\ \ref{tdrf3.10-1}. $\;$Fix $Y\in \cy$.
One can  choose $\eta\in C_c(\ce)$  with $0\leq \eta\leq 1$ and $\eta=1$ on $Y$ by  Lemma~\ref{l2}.\ref{l2-2}.
Then if $\psi\in B(\ce)$ one has $\eta\,\psi\in B(\ce_M)$.
Moreover,
\[
(\ce_M)_{m,Y}(\eta\,\psi)=\sup_{\xi\in \cc_Y(\ce_M)}(\ce_M)_\xi(\eta\,\psi)
=\sup_{\xi\in \cc_Y(\ce)}\ce_\xi(\eta\,\psi)=\sup_{\xi\in \cc_Y(\ce)}\ce_\xi(\psi)=\ce_{m,Y}(\psi)
\]
where the third step uses locality.
Next replacing $\eta$ by a net $\eta_Z\in \cc_Z(\ce)$ with  $\eta_Z=1$ on $Y$ and $\lim_{Z\in\cy}\eta_Z=1$ pointwise
one deduces that 
\[
\liminf_{Z\in\cy}(\ce_M)_{m,Y}(\eta_Z\psi)=\ce_{m,Y}(\psi)<\infty
\]
for all $\psi\in B (\ce)$.
But $\lim_{Z\in\cy}\|\eta_Z\psi-\psi\|_2=0$.
Therefore  $\psi$ is in the domain of the relaxation 
$(\ce_M)_{m,Y;\,0}$ and 
\[
(\ce_M)_{m,Y;\,0}(\psi)=\liminf_{\varphi\to\psi}(\ce_M)_{m,Y}(\varphi)
\leq \liminf_{Z\in\cy}(\ce_M)_{m,Y}(\eta_Z\psi)=\ce_{m,Y}(\psi)
\]
for all $\psi\in B (\ce)$.
But $\ce_{m,Y}\leq\ce$ and $B (\ce)$ is a core of $\ce$.
Therefore the inequality  $(\ce_M)_{m,Y;\,0}(\psi)\leq \ce_{m,Y}(\psi)$ extends to all $\psi\in D(\ce)=D(\ce_{m,Y})$ by continuity. 
Consequently 
 $(\ce_M)_{m,Y;\,0}\leq \ce_{m,Y;\,0}$. 
 Conversely $\ce\leq \ce_M$ and $\ce_\xi\leq (\ce_M)_\xi$ for all $\xi\in \cc_Y(\ce)=\cc_Y(\ce_M)$.
 Therefore $\ce_{m,Y}\leq (\ce_M)_{m,Y}$ and correspondingly
 $\ce_{m,Y;\,0}\leq (\ce_M)_{m,Y;\,0}$.
 Combining these conclusions one has
 \begin{equation}
 \ce_{m,Y;\,0}=(\ce_M)_{m,Y;\,0}
 \label{edrf3.11}
 \end{equation}
 for all $Y\in \cy$.
 Then taking the limit over $Y\in\cy$ one obtains the first statement of the theorem, $\ce_m=(\ce_M)_m$.
 
  \begin{remarkn}\label{rdrf3.11} The identity (\ref{edrf3.11})  establishes that $\ce_M$ is a core of $\ce_{m,Y}$.
 Specifically  (\ref{edrf3.11})  states that the relaxation of  $\ce_{m,Y}$ restricted to $D(\ce_M)$ is equal to  $\ce_{m,Y;\,0}$.
 Then each core of $\ce_M$ is a core of $\ce_{m,Y}$.
 \end{remarkn}  
 
 \smallskip

\noindent \ref{tdrf3.10-2} and \ref{tdrf3.10-3}. $\;$First note that $\ce_{m,Y}(\varphi)=\ce(\varphi)=\ce_M(\varphi)$ for all $\varphi\in \cc_Y(\ce)$.
But it is not clear that this identity extends to the relaxation $\ce_{m,Y;0}$
It does, however, for $\varphi$ with support strictly in the interior of $Y$.

Let $\varphi\in \cc_Z(\ce)$ where $Z\in\cy$ is such that $\overline Z\subset Y$.
Then $\ce_{m,Y;0}(\varphi)=\liminf_{\psi\to\varphi}\ce_{m,Y}(\psi)$.
In particular $\psi=\{\psi_n\}_{n\geq1}$ is a sequence of $\psi_n\in D(\ce)$ such that $\|\psi_n-\varphi\|_2\to0$
and $\liminf_{n\to\infty}\ce_{m,Y}(\psi_n)<\infty$.
Now one can choose $\eta\in \cc_Y(\ce)$ with $0\leq \eta\leq 1$ and $\eta=1$ on $ Z$.
Then $\eta\,\psi_n\in \cc_Y(\ce)$.
Moreover, since $\eta\,\varphi=\varphi$ it follows that $\|\eta\,\psi_n-\varphi\|_2\to0$.
But one again deduces by locality that 
\[
\ce_{m,Y}(\eta\,\psi_n)=\sup_{\xi\in \cc_Y(\ce)}\ce_\xi(\eta\,\psi_n)=\sup_{\xi\in \cc_Y(\ce)}\ce_\xi(\psi_n)=\ce_{m,Y}(\psi_n)
\;.
\]
Therefore $\liminf_{n\to\infty}\ce_{m,Y}(\eta\,\psi_n)=\liminf_{n\to\infty}\ce_{m,Y}(\psi_n)<\infty$.
One then has
\begin{eqnarray*}
\ce_{m,Y;0}(\varphi)=\liminf_{\psi\to\varphi}\ce_{m,Y}(\psi)&=&\liminf_{\psi\to\varphi}\ce_{m,Y}(\eta\,\psi)\\[5pt]
&=&\liminf_{\psi\to\varphi}\ce_M(\eta\,\psi)\geq
\liminf_{\psi\to\varphi}\ce_M(\psi)=\ce_M(\varphi)
\end{eqnarray*}
where the third equality follows because the $\eta\,\psi_n\in \cc_Y(\ce)$.
Consequently $\ce_m(\varphi)\geq\ce_M(\varphi)$ for all $\varphi\in \cc_Z(\ce)$ and all $Z\in\cy$.
But $\ce_m\leq \ce_M$ so this latter relation must be an equality.
Therefore   $\ce_m(\varphi)=\ce_M(\varphi)$ for all $\varphi\in C_c(\ce_M)$ and by closure for  all $\varphi\in D(\ce_M)$.
Hence $\ce_m\supseteq \ce_M$ and $(\ce_m)_M=\ce_M$.

 \smallskip
 
 \noindent \ref{tdrf3.10-4}. $\;$The proof of this statement is somewhat longer but it is the key to proving the remaining statements.
It depends on  the order properties of the Markovian semigroups  associated with the various Dirichlet forms.

 First we introduce  forms $\{\ce_{M,Y}\}_{Y\in\cy}$ analogous to the forms defined in the construction  $ \ce_m$.
For each   $Y\in\cy$ let  $C_Y(\ce)=D(\ce)\cap C_c(Y)$.
Then  define $\ce_{M,Y}$ as  the closure of the restriction of $\ce$ to $C_Y(\ce)$.
Note that $C_Y(\ce)$ is dense in $L_2(Y)$ by Lemma~\ref{l2}.\ref{l2-1}.
Therefore $\ce_{M,Y}$ is a densely defined regular Dirichlet form  on $L_2(Y)$.

Secondly let $(H_{M,Y}, S^{M,Y})$ denote the operator and semigroup on $L_2(Y)$ canonically associated with 
$\ce_{M,Y}$ on  $L_2(Y)$.
Similarly let $(H_{M}, S^{M})$ and $(H, S)$ denote the operator and semigroup on $L_2(X)$ associated with 
 $\ce_M$ and $\ce$, respectively.
Next we adopt the convention that each bounded operator $B$ on the subspace $L_2(Y)$ can be extended to a bounded operator
on $L_2(X)$, still denoted by $B$, through the definition $\varphi\in L_2(X)\mapsto B(\one_Y\varphi)\in L_2(Y)$.
In particular the resolvents $(\lambda I+H_{M,Y})^{-1}$ and the semigroups $S^{M,Y}$  extend from $L_2(Y)$ to  $L_2(X)$.
The extended semigroups $S^{M,Y}$ are  strongly continuous  on $L_2(X)$ but  $S^{M,Y}_t\to\one_Y$ as $t\to0$.
In the sequel we use the observation that Propositions~\ref{pdrf2.4} and \ref{pdrf2.5} are  valid for the semigroups $S^{M,Y}$
extended to $L_2(X)$  (see, for example,  \cite{bSim5} or \cite{Ouh5}, Section~2.6).

\smallskip

Next we consider the order relations for the extended  semigroups  on $L_2(X)$.

\begin{prop}\label{pdrf3.11}
If $\ce$ is a  local, semi-regular,   Dirichlet form and $Y_1, Y_2\in\cy $ with $Y_1\subseteq Y_2$ then
\begin{equation}
0\leq S^{M,Y_1}_t\varphi\leq S^{M,Y_2}_t\varphi\leq S^M_t\varphi
\label{emf1}
\end{equation}
for all $\varphi\in L_2(X)_+$ and all $t>0$.
Moreover, the monotonically increasing net of semigroups $\{S^{M,Y}\}_{Y\in\cy}$  converges strongly to $S^M$.

If, in addition, $\ce$ is inner regular then 
\begin{equation}
0\leq  S^M_t\varphi\leq S_t\varphi
\label{emf11}
\end{equation}
for all $\varphi\in L_2(X)_+$ and all $t>0$.
\end{prop}
\proof\
The  ordering of the $S^{M,Y}$  is proved by a refinement of the argument of \cite{ER30}.
In particular it depends on the following  variation of Lemma~2.2 in this reference.

\begin{lemma}\label{lmf1}
Assume $\ce$ is a  local, semi-regular,   Dirichlet form.
If $\varphi\in D(\ce_{M,Y_1})$, $\psi\in D(\ce_{M,Y_2})_+$  where $Y_1\subseteq Y_2$ and if
\begin{equation}
(\chi,\varphi)+\ce_{M,Y_1}(\chi,\varphi)\leq (\chi,\psi)+\ce_{M,Y_2}(\chi,\psi)
\label{emf2}
\end{equation}
for all $\chi\in D(\ce_{M,Y_1})_+$ then $\varphi\leq \psi$ .
\end{lemma}
\proof\
Since $\varphi\in D(\ce_{M,Y_1})$ and $\psi\in D(\ce_{M,Y_2}) $ there exist sequences $\varphi_n \in C_{Y_1}(\ce)$ and $\psi_m\in C_{Y_2}(\ce)$ with 
$\lim \|\varphi_n - \varphi\|_{D(\ce)} = 0$ and  $\lim \|\psi_m - \psi\|_{D(\ce)} = 0$.
But $\psi\geq 0$ and  since the modulus map $\psi_m\mapsto |\psi_m|$ is continuous in the $D(\ce)$-graph norm
it follows that  $\lim \||\psi_m| - \psi\|_{D(\ce)} = 0$.
Therefore one may assume that the $\psi_m\in C_{Y_2}(\ce)_+$.
Then, however, $\supp (\varphi_n - \psi_m)_+ \subseteq \supp \varphi_n \subset Y_1$
since $\psi _m\geq 0$.
So $(\varphi_n - \psi_m)_+ \in D(\ce_{M,Y_1})$ for all $n,m$.
Moreover, $\lim \|(\varphi_n - \psi_m)_+ - (\varphi - \psi)_+\|_{D(\ce)}=0$.
Hence $(\varphi - \psi)_+ \in D(\ce_{M,Y_1} )$.

Secondly, set $\chi=(\varphi-\psi)_+$ in  (\ref{emf2}).
Then one deduces that 
\[
\|(\varphi-\psi)_+\|_2^2
=((\varphi-\psi)_+, \varphi-\psi)
\leq -\ce((\varphi-\psi)_+, \varphi-\psi)
=-\ce((\varphi-\psi)_+)
\leq 0
\;\;\;,
\]
where we used locality of $\ce$ in the last equality.
Hence $(\varphi-\psi)_+ = 0$ or, equivalently, $\varphi \leq \psi$.\hfill$\Box$

\bigskip
The first ordering property of Proposition~\ref{pdrf3.11} now follows as in \cite{ER30}.

\smallskip

Fix  $\eta\in L_2(X)_+$. 
Then  set $\varphi=(I+H_{M,Y_1})^{-1}\eta $ and $\psi=(I+H_{M,Y_2})^{-1}\eta$.
Hence   $\varphi\in D(\ce_{M,Y_1})$ and $\psi\in D(\ce_{M,Y_2})_+$.
Moreover, (\ref{emf2}) is satisfied. 
Therefore $(I+H_{M,Y_1})^{-1}\eta\leq (I+H_{M,Y_2})^{-1}\eta$.
By rescaling and iterating one concludes that $(I+t\,H_{M,Y_1}/n)^{-n}\eta\leq (I+t\,H_{M,Y_2}/n)^{-n}\eta$ for all $t>0$ and $n$.
Therefore $S^{M,Y_1}_t\eta\leq S^{M,Y_2}_t\eta$ for all $t>0$ by the Trotter product formula.

\smallskip

The second order property  of the proposition follows from the convergence statement  which is established as follows.

The net of submarkovian semigroups $\{S^{M,Y}\}_{Y\in\cy}$ is  monotonically increasing by the first statement.
It then follows as a corollary of an old result of Vigier (see \cite{RN}, page~261) that the net converges strongly on $L_2(X)$ to a submarkovian semigroup~$T$.
In fact the $S^{M,Y}$ converge 
strongly to $T$ on each of the spaces $L_p(X)$, $p\in[1,\infty\rangle$.
The $L_p$-convergence follows from results of Karlin \cite{Karl} and Krasnoselski \cite{Kra}
(see \cite{KisR2} Propositions A3 and A4).
It remains to identify $T$ and $S^M$.

Let $H$ denote the generator of $T$ and $\cf$ the Dirichlet form
corresponding to $H$.
The net of forms $\{\ce_{M,Y}\}_{Y\in\cy}$ is  monotonically decreasing.
Hence it follows from Proposition~\ref{pdrf2.5} that 
  $\cf$ is the relaxation of the 
limit form $\ce_{M,X}$ defined by $D(\ce_{M,X})=\bigcup_{Y\in\cy} D(\ce_{M,Y})$ and 
\[
\ce_{M,X}(\varphi)=\lim_{Y\in\cy}\ce_{M,Y}(\varphi)
\]
for all $\varphi\in D(\ce_{M,X})$.
But $D(\ce_{M,X})\supseteq C_c(\ce)$ and 
 $\ce_{M,X}(\varphi)=\ce_M(\varphi)$ for all $\varphi\in C_c(\ce)$
 and therefore by continuity for all $\varphi\in D(\ce_M)$.
Hence  $\cf\leq \ce_M$.

Secondly, since the net of semigroups $\{S^{M,Y}\}_{Y\in\cy}$ converges strongly to $T$ the corresponding  net of generators  $Y\in\cy\mapsto H_{M,Y}$ converges in the strong resolvent sense to $H$.
But $\varphi\in L_2(X)\mapsto (I+H_{M,Y})^{-1}\one_Y\varphi\in D(H_{M,Y})\subseteq D(\ce_{M,Y})\subseteq D(\ce_M)$ and
\[
\|(I+H_{M,Y})^{-1}\one_Y\varphi\|_{D(\ce_M)}^2=(\one_Y\varphi,(I+H_{M,Y})^{-1}\one_Y\varphi)\leq \|\varphi\|_2^2
\;.
\]
Therefore  it follows from the Banach--Alaoglu theorem (see, for example, \cite{Ouh5}, Lemma~1.32, or \cite{MR}, Lemma~I.2.12) that $(I+H)^{-1}\varphi\in D(\ce_M)$ and 
\begin{eqnarray*}
\|(I+H)^{-1}\varphi\|_{D(\ce_M)}^2&\leq&\liminf_{Y\in\cy}\|(I+H_{M,Y})^{-1}\one_Y\varphi\|_{D(\ce_M)}^2\\[5pt]
&=&(\varphi,(I+H)^{-1}\varphi)
=\|(I+H)^{-1}\varphi\|_{D(\cf)}^2
\;.
\end{eqnarray*}
But this immediately implies that $\ce_M(\psi)\leq \cf(\psi)$ for all $\psi\in D(H)$.
Since $D(H)$ is a core of $\cf$ one then deduces that $\ce_M\leq \cf$.
As we have already established that  $\cf\leq \ce_M$ it follows that the two forms are equal.
Therefore $T=S^M$ and one concludes  that  $\{S^{M,Y}\}_{Y\in\cy}$ converges strongly to $S^M$.
Since the net of semigroups is monotonically increasing one immediately deduces that $0\leq S^{M,Y}\varphi\leq S^M_t\varphi$
for all $Y\in\cy$, $\varphi\in L_2(X)_+$ and $t\geq0$.

\smallskip

It remains to establish the domination property (\ref{emf11}) for inner regular $\ce$.
But $\ce_M$ is the closure with respect to the graph norm $\|\,\cdot\,\|_{D(\ce)}$ of $\ce$
restricted to $C_c(\ce)$.
Then by inner regularity it is the closure of $\ce$ restricted to $B_c(\ce)$ or, equivalently, to the subspace $D_c(\ce)$ of $D(\ce)$
spanned by the functions with compact support.
Now if $\varphi\in D(\ce_M)$,  $\psi\in D(\ce)_+$ and 
\[
(\chi,\varphi)+\ce_{M}(\chi,\varphi)\leq (\chi,\psi)+\ce(\chi,\psi)
\]
for all $\chi\in D(\ce_M)_+$ then $\varphi\leq \psi$ by Lemma~2.2 of \cite{ER30}.
But then (\ref{emf11}) follows exactly as in the proof of Proposition~2.1 in \cite{ER30}. 
\hfill$\Box$

\bigskip

The next proposition gives an ordering of  the semigroups $S^{m, Y;\,0}$ and $S^{M, Y}$  corresponding to the forms   $\ce_{m,Y;\,0}$  and $\ce_{M, Y}$, respectively.
The proof is a variation of the argument used to prove Proposition~3.4  in \cite{ER30}.

\begin{prop}\label{pdrf3.13}
If $\ce$ is a local,    inner regular, Dirichlet form and $Y\in \cy$   then
\begin{equation}
0\leq S^{M,Y}_t\varphi\leq  S^{m,Y;\,0}_t\varphi
\label{edfr3.1}
\end{equation}
for all $\varphi\in L_2(X)_+$ and all $t>0$.
\end{prop}
\proof\
Define the family of forms $\{\ce_{m, Y;\,\varepsilon}\}_{\varepsilon>0}$  on the common domain $D(\ce)$ by 
$\ce_{m,Y;\,\varepsilon}=\ce_{m,Y}+\varepsilon\,\ce$.
Since $\ce_{m,Y}\leq  \ce$ it follows that  $\varepsilon \,\ce\leq \ce_{m,Y;\,\varepsilon}\leq (1+\varepsilon)\,\ce$.
Therefore the forms $\ce_{m,Y;\,\varepsilon}$ are closed on  $D(\ce)$.
Moreover, the $\ce_{m,Y;\,\varepsilon}$  are local because the $\ce_{m,Y}$ inherit locality from $\ce$.  
Then since $B_c(\ce_{m,Y;\,\varepsilon})=B_c(\ce)$ the  $\ce_{m,Y;\,\varepsilon}$  are inner regular Dirichlet forms.  
But 
\[
C_Y(\ce_{m,Y;\,\varepsilon})=D(\ce_{m,Y;\,\varepsilon})\cap C_c(Y)=D(\ce)\cap C_c(Y)=C_Y(\ce)\;.
\]
Hence $(\ce_{m,Y;\,\varepsilon})_{M, Y}=(1+\varepsilon)\,\ce_{M, Y}$.
Next, for brevity,  let $T$ denote the submarkovian semigroup $S^{m,Y;\,\varepsilon}$   associated with the Dirichlet form $\ce_{m,Y;\,\varepsilon}$.
Then the foregoing identification implies that $S^{M,Y}_{(1+\varepsilon)t} =T_t^{M,Y}$ for all $t\geq0$.
 Since the  $\ce_{m,Y;\,\varepsilon}$ are both local and  inner regular $0\leq T^{M,Y}_t\varphi\leq T^M_t\varphi\leq T_t\varphi$ for all $\varphi\in L_2(X)_+$ and $t>0$ 
 by the last statement of Proposition~\ref{pdrf3.11} applied  with $S$ replaced by  $T$.
Combining these conclusions and substituting $T_t=S_t^{m,Y;\,\varepsilon}$ one deduces  that 
\begin{equation}
0\leq S^{M,Y}_{(1+\varepsilon)t}\varphi \leq S_t^{m,Y;\,\varepsilon}\varphi
\label{edrf3.2}
\end{equation}
 for all $\varphi\in L_2(X)_+$ and $t>0$.

 The forms $\ce_{m,Y;\, \varepsilon}$ decrease monotonically as $\varepsilon\to0$.
 Therefore it follows from Proposition~\ref{pdrf2.5} that 
 the corresponding positive self-adjoint operators $H_{m, Y;\,\varepsilon}$  converge in the strong
resolvent sense to the  positive self-adjoint operator $H_Y$  associated with the relaxation of the form $h_Y(\varphi)=\lim_{\varepsilon\to0}\ce_{m,Y;\,\varepsilon}(\varphi)$ on the  domain 
$D(h_Y)=\bigcup_{\varepsilon>0}D(\ce_{m,Y;\,\varepsilon})$.
But $D(h_Y)=D(\ce)$ and $h_Y(\varphi)=\ce_{m,Y}(\varphi)$ for $\varphi\in D(\ce)$.
Therefore $H_Y$ is the operator corresponding to the relaxation $\ce_{m,Y;\,0}$ of $\ce_{m,Y}$ since $D(\ce)$ is a core of $\ce_{m,Y}$ by Remark~\ref{rdrf3.11}.
Consequently $S^{m,Y;\,\varepsilon}_t\varphi\to S^{m,Y;\,0}_t\varphi$ as $\varepsilon\to0$.
Since $S^{M,Y}_{(1+\varepsilon)t}\varphi\to S^{M,Y}_{t}\varphi$ as $\varepsilon\to0$ the assertion (\ref{edfr3.1})  follows immediately from (\ref{edrf3.2}).
\hfill$\Box$

 \bigskip

Now we can establish the key order relation for $S^M$ and $S^m$.

\begin{prop}\label{pdrf3.14}
If $\ce$ is a local,  inner regular, Dirichlet form and $S^M$, $S^m$ are the submarkovian semigroups
corresponding to the forms $\ce_M$, $\ce_m$ then
\begin{equation}
0\leq S^M_t\varphi\leq S^m_t\varphi
\label{edfr3.3}
\end{equation}
for all $\varphi\in L_2(X)_+$ and all $t>0$.
Consequently, $D(\ce_M)$ is an order ideal of $D(\ce_m)$ and 
$B(\ce_M)$ is an algebraic ideal of $B(\ce_m)$.
Therefore 
$B_c(\ce)=B_c(\ce_M)=B_c(\ce_m)$ and $C_c(\ce)=C_c(\ce_M)=C_c(\ce_m)$.
\end{prop}
\proof\
It follows from Proposition~\ref{pdrf3.11} that the net of semigroups
 $\{S^{M,Y}\}_{Y\in\cy}$ converges  strongly to $S^M$.
Next the generator of the semigroup $S^m$ is the operator associated with the Dirichlet form $\ce_m$ defined by (\ref{edfr3.0}).
Therefore $\ce_m$ is the supremum   of the monotonically increasing net $\{ \ce_{m,Y;\,0}\}_{Y\in\cy}$ of Dirichlet forms $\ce_{m,Y;\,0}$.
Hence the corresponding semigroups $S^{m,Y;\,0}$ converge strongly to the semigroup $S^m$ by Proposition~\ref{pdrf2.4}.
Then the semigroup ordering follows from  (\ref{edfr3.1}).
Explicitly one has
\[
0\leq S^M_t\varphi=\lim_{Y\in\cy}S^{M,Y}\varphi\leq \lim_{Y\in\cy}S^{m,Y;\,0}\varphi=S^m_t\varphi
\]
for all $\varphi\in L_2(X)_+$ and all $t>0$.

The assertions that  $D(\ce_M)$ is an order ideal of $D(\ce_m)$  and $B(\ce_M)$ is an algebraic ideal of $B(\ce_m)$ now follow from 
Proposition~\ref{tdrf2} with $\ce$ replaced by $\ce_M$, $\cf$ by $\ce_m$, $S_t$ by $S^M_t$ and $T_t$ by $S^m_t$ and noting that $\ce_m\supseteq \ce_M$.
Then  Proposition~\ref{tdrf2} establishes that the order property  $0\leq S^M_t\varphi\leq S^m_t\varphi$ for all $\varphi\geq 0$ and $t>0$  is equivalent to the ideal properties .
 
 The  last statement of the proposition follows  because  $B(\ce_M)$ is an algebraic ideal of $B(\ce_m)$.
 In particular $B_c(\ce_M)\,B_c(\ce_m)\subseteq B_c(\ce_M)\,B(\ce_m) \subseteq B_c(\ce_M)$.
Now if $\varphi\in B_c(\ce_m)$ it follows from  Lemma~\ref{l2}.\ref{l2-2}   that one can choose $\eta\in C_c(\ce_M)$ such that $0\leq \eta\leq \one_X$ 
and $\eta=1$ on the support of $\varphi$. Therefore $\eta\,\varphi=\varphi$.\noindent
Hence $\varphi\in C_c(\ce_M) B_c(\ce_m)\subseteq  B(\ce_M)$.
Thus $B_c(\ce_m)\subseteq B_c(\ce_M)$.
But $\ce_m\leq \ce\leq \ce_M$ and consequently $B_c(\ce_m)\supseteq B_c(\ce)\supseteq B_c(\ce_M)$.
Therefore $B_c(\ce)=B_c(\ce_M)=B_c(\ce_m)$.
Finally $C_c(\ce)=B_c(\ce)\cap C(X)$ etc.\ so one also deduces that $C_c(\ce)=C_c(\ce_M)=C_c(\ce_m)$.
\hfill$\Box$

\bigskip

\noindent{\bf Proof of Theorem~\ref{tdrf3.10} continued}$\;$    
Statements~\ref{tdrf3.10-4}   and \ref{tdrf3.10-41}   of Theorem~\ref{tdrf3.10} are established by the last statement of Proposition~\ref{pdrf3.14}.
It remains to prove Statement~\ref{tdrf3.10-5}. 

\smallskip

\noindent \ref{tdrf3.10-5}. $\;$ First $\ce_M$ is regular and 
$C_c(\ce_M)=C_c(\ce)$.
Secondly, the inner regularity of $\ce_m$ follows 
from a slightly more general result.

\begin{lemma}\label{inreg} Let  $\ce$ be  an inner regular Dirichlet form and 
$\cf$  a  Dirichlet form extension of $\ce_M$.
Assume $B_c(\cf)=B_c(\ce_M)$.

It follows that $\cf$ is inner regular and $\cf_M=\ce_M$.
\end{lemma} 
\proof\
Since $B_c(\cf)=B_c(\ce_M)$ one also has $C_c(\cf)=C_c(\ce_M)$.
But 
$\ce_M=\overline{\cf|_{C_c(\ce_M)}}$ because $\cf$ is an extension of $\ce_M$.
Then
\[
\ce_M=\overline{\cf|_{C_c(\ce_M)}}=\overline{\cf|_{C_c(\cf)}}
\subseteq \overline{\cf|_{B_c(\cf)}}=\overline{\cf|_{B_c(\ce_M)}}
=\overline{\ce_M|_{B_c(\ce_M)}}=\overline{\ce_M|_{C_c(\ce_M)}}=\ce_M
\]
where the penultimate step uses the  regularity of $\ce_M$.
Therefore $\overline{\cf|_{C_c(\cf)}}=\overline{\cf|_{B_c(\cf)}}=\cf_M$ and $\cf$ is inner regular.
Moreover, $\cf_M=\overline{\cf|_{C_c(\cf)}}=\ce_M$.
\hfill$\Box$

\bigskip
The inner regularity of $\ce_m$ now  follows from setting $\cf=\ce_m$ in Lemma~\ref{inreg}  and noting that the assumptions of the lemma
are satisfied by Statements~\ref{tdrf3.10-3} and \ref{tdrf3.10-4} 
of the theorem.
This completes the proof of Theorem~\ref{tdrf3.10}.\hfill$\Box$

\bigskip

Since  $\ce_M$ is a restriction of $\ce$ it follows that $\ce_M$ inherits the locality property from $\ce$. 
It is,  however, unclear if $\ce_m$ is  local. 
This does follow if the forms $\ce_{m,Y}$ are closable for all $Y\in\cy$ which in turn follows if the
 truncated forms $\ce_\xi$ are closable for all $\xi\in \cc_Y(\ce)$ and $Y\in\cy$.

\begin{prop}\label{pdrf4.1}
Let $\ce$ be a local $($resp.\ strongly local\,$)$, semi-regular, Dirichlet form.
If the forms $\ce_{m,Y}$ are closable for all $Y\in\cy$
then $\ce_m$ is local $($resp.\ strongly local\,$)$.
\end{prop}
\proof\
Since $\ce_{m,Y}$ is closable the relaxation $\ce_{m,Y;\,0}$ is  automatically equal to its closure $\overline{\ce}_{m,Y}$.
Then $\ce_m$ is the monotonic limit of the increasing net   $Y\in\cy\mapsto \overline{\ce}_{m,Y}$.
Thus 
\begin{equation}
\ce_m(\varphi)=\lim_{Y\in\cy}\overline{\ce}_{m,Y}(\varphi)=\sup_{Y\in\cy}\overline{\ce}_{m,Y}(\varphi)
\label{edrf4.111}
\end{equation}
for all $\varphi\in D(\ce_m)=\{\varphi\in\bigcap_{Y\in \cy}D(\overline{\ce}_{m,Y}):\sup_{Y\in\cy}\overline{\ce}_{m,Y}(\varphi)<\infty\}$.

Secondly, if  $\ce$ is local  the truncated functions $\ce_\xi$ are local on $B(\ce)$.
Therefore $\ce_\xi(|\varphi|)=\ce_\xi(\varphi)$ for all $\xi\in \cc_Y(\ce)$ and all $\varphi\in B(\ce)$.
Hence
\[
\ce_{m,Y}(|\varphi|)=\lim_{\xi\in \cc_Y(\ce)}\ce_\xi(|\varphi|)=\lim_{\xi\in \cc_Y(\ce)}\ce_\xi(\varphi)=\ce_{m,Y}(\varphi)
\]
for all $\varphi\in B(\ce)$.
It then follows from Proposition~\ref{pdrf2.3}, applied with $\ce={\overline\ce}_{m,Y} $ and $D=B(\ce)$, that  $\overline{\ce}_{m,Y}$ is local.
Then  $\ce_m$ is local  by a similar argument.

The proof is analogous if  $\ce$ is strongly local.
The locality criterion $\ce(|\varphi|)=\ce(\varphi)$ of Proposition~\ref{pdrf2.3} is replaced 
by the corresponding criterion $ \ce(|\varphi+1|-1)=\ce(\varphi)$  for strong locality.
\hfill$\Box$

\begin{remarkn}\label{rdrf3.1} 
The foregoing argument only requires locality, or strong locality, of the
approximants $\{\ce_{m,Y}\}_{Y\in \cy}$.
This of course follows from locality, or strong locality,  of $\ce$ although the latter property is not necessary. 
It is quite possible that $\ce$ is local but the $\ce_{m,Y}$ are strongly local.
In particular $\ce_m$ can be strongly local even if $\ce$ is only local (see the discussion of Robin boundary conditions
at the end of Section~\ref{S4}).
\end{remarkn}

The assumption that the forms $\ce_{m,Y}$ are closable for all $Y\in\cy$ is not very satisfactory 
although it is satisfied for a large class of locally strongly elliptic operators.
In fact the  $\ce_{m,Y}$ are often closed.
For example, consider the Laplacian defined on the open subset $\Omega$ of $\Ri^d$ with domain
$C_c^\infty(\Omega)$.
Then $\ce_{m,Y}(\varphi)=\|\nabla\varphi\|_2^2$ with domain consisting of the $\varphi\in L_2(\Omega)$ whose restriction to $Y$ is in  $W^{1,2}(Y)$.
It follows that  $\ce_{m,Y}$ is closable but if the boundary of $Y$ is smooth, for example Lipschitz, then  $\ce_{m,Y}$ is closed.

\bigskip

Finally we give two examples of elliptic operators in one-dimension which illustrate the conclusions of Theorem~\ref{tdrf3.10}.
These examples are analyzed in detail in \cite{RSi3}.

\begin{exam}\label{ex1d3.1}
Define the positive, symmetric, operator $H$ on $L_2(\Ri_+)$ by  $H\varphi=-(c\,\varphi')'$ where $c\in W^{1,\infty}_{\rm loc}(\Ri_+)$ is strictly positive  
and $\varphi\in D(H)=C_c^\infty(\Ri_+)$.
Let  $\ce_0(\varphi)=\int^\infty_0 c\,(\varphi')^2$ with $D(\ce_0)=C_c^\infty(\Ri_+)$.
Set $\nu(x)=\int_x^1 c^{-1}$.

\smallskip
 
 1. If $\nu\in L_\infty(0,1)$ then $H$ has a one-parameter family of submarkovian extensions
 $H^{(\alpha)}$, where  $\alpha\in[0,\infty]$, with  corresponding Dirichlet forms 
 $\ce^{(\alpha)}$.
 If $\alpha\in[0,\infty\rangle$ then 
 $\ce^{(\alpha)}(\varphi) =\ce^{(0)}(\varphi)+\alpha\,|\varphi(0)|^2$ 
  where $\ce^{(0)}$ is the extension of $\ce_0$ to the domain $D(\overline{\ce_0})+\spann\sigma_+$
  with $\sigma_+\in C_c^\infty(\Ri_+)$, $0\leq \sigma_+\leq 1$, $\sigma_+(x)=1$ if $ x\in\langle0,1\rangle$ and $\sigma_+( x)=0$
  if $ x\geq 2$ and where $D(\ce^{(\alpha)})=D(\ce^{(0)})$.
 The associated operators satisfy the boundary conditions $(c\,\varphi')(0)=\alpha\,\varphi(0)$.
 In addition $\ce^{(\infty)}=\overline{\ce_0}$ and the corresponding boundary condition  is  $\varphi(0)=0$.
 The family of forms $\ce^{(\alpha)}$ is monotonically increasing with $\alpha$
and  $\ce^{(\infty)}$ formally corresponds
 to the limit $\alpha\to\infty$ of the $\ce^{(\alpha)}$.
  These statements are contained in \cite{RSi3}, Theorem~2.4.
 
  The forms $\ce^{(\alpha)}$ are all inner regular and $\ce^{(\infty)}$ is regular.
One can calculate explicitly the maximal and minimal forms $(\ce^{(\alpha)})_M$ and  $(\ce^{(\alpha)})_m$.
 First it follows that $C_c(\ce^{(\alpha)})=C_c(\overline{\ce_0})$ for all $\alpha\in[0,\infty]$.
Therefore   $(\ce^{(\alpha)})_M=\overline{\ce_0}=\ce^{(\infty)}$ for all $\alpha\in[0,\infty]$.
Secondly, $\ce^{(\alpha)}_\xi(\varphi)=\int_Y\xi \,c\,(\varphi')^2$ for $\xi\in C_c(\ce^{(\alpha)})$ with $\supp \xi\subseteq Y$ 
and $\varphi\in B(\ce^{(\alpha)})$.
Hence $\ce^{(\alpha)}_{m,Y}(\varphi)=\int_Yc\,(\varphi')^2$ for all $\varphi\in B(\ce^{(\alpha)})$.
But if $\alpha\in[0,\infty\rangle$ then $B(\ce^{(\alpha)})=B(\ce^{(0)})$.
Alternatively  $\ce^{(\infty)}_{m,Y}$ extends to $B(\ce^{(0)})$ by continuity.
Then by definition $(\ce^{(\alpha)})_m=\ce^{(0)}$ for all $\alpha\in[0,\infty]$.
Thus the maximal and minimal forms are independent of the choice of $\alpha$ and are identified 
with the maximal and minimal forms in the family  $\ce^{(\alpha)}$.
The minimal form $\ce^{(0)}$ is an extension of the maximal form $\ce^{(\infty)}$ but not of the intermediate
forms $\ce^{(\alpha)}$ with $\alpha\in\langle0,\infty\rangle$.

 \smallskip
 
 2. If $\nu\not\in L_\infty(0,1)$ then $H$ has a unique submarkovian extension determined by the form $\overline{\ce_0}$.
 The associated self-adjoint operator satisfies the boundary condition $(c\,\varphi')(0)=0$.
 In this case $\ce_M=\ce_m=\overline{\ce_0}$.

  \end{exam}
  
  The diffusion process described by the first case of the example  is determined by its behaviour at the boundary point, i.e.\ at the origin.
  This behaviour is dependent on the boundary value of the diffusion coefficient $c$.
  In the second case the coefficient is sufficiently degenerate at the origin that the diffusion fails to reach the boundary.
  Therefore there is no ambiguity and the diffusion is uniquely determined.
  A different phenomenon can occur  if the diffusion is degenerate at an interior point.
  This is illustrated by the second example which also clarifies the significance of inner regularity.  
    \begin{exam}\label{ex1d3.2}
Define the positive, symmetric, operator $H$ on $L_2(\Ri)$ by  $H\varphi=-(c\,\varphi')'$ where $\varphi\in D(H)=C_c^\infty(\Ri)$
and $c\in W^{1,\infty}_{\rm loc}(\Ri)$ is strictly positive  on $\Ri\backslash\{0\}$   but $c(0)=0$.
Let  $\ce_0(\varphi)=\int^\infty_{-\infty} c\,(\varphi')^2$ with $D(\ce_0)=C_c^\infty(\Ri)$.
Set $\nu_+(x)=\int_x^1 c^{-1}$ and $\nu_-=\int_{-1}^{-x}c^{-1}$.

\smallskip

 1. If $\nu_+\vee \nu_-\in L_\infty(0,1)$ then $H$ has a one-parameter family of submarkovian extensions
 $H^{(\alpha)}$, where  $\alpha\in[0,\infty]$.
 If $\alpha\in[0,\infty\rangle$ then these extensions are determined by the  Dirichlet forms 
 \[
 \ce^{(\alpha)}(\varphi)=\ce^{(0)}(\varphi)+\alpha\,|\varphi(0_+)-\varphi(0_-)|^2
 \]
 where $\ce^{(0)}$ is the extension of $\ce_0$ to the domain $D(\overline{\ce_0})+\spann(\sigma_+-\sigma_-)$
  and $D(\ce^{(\alpha)})=D(\ce^{(0)})$.
Here $\sigma_+$ is the function on $\Ri_+$ defined in Example~\ref{ex1d3.1} and $\sigma_-$ is defined on $\Ri_-$ by
 $\sigma_-(x)=\sigma_+(-x)$.   
 The domain of the corresponding self-adjoint  operator is characterized by the continuity condition $(c\,\varphi')(0_+)-(c\,\varphi')(0_-)=\alpha\,(\varphi(0_+)-\varphi(0_-))$.
 In addition  $\ce^{(\infty)}=\overline{\ce_0}$ which again corresponds to the limit $\alpha\to\infty$.
 (These statements are established in \cite{RSi3}, Theorem~1.1.)
 One can again calculate the maximal and minimal forms.

First  $C_c(\ce^{(\alpha)})=C_c(\overline{\ce_0})$ for all $\alpha\in[0,\infty]$.
(Note that here it is essential that we are considering the subspace of 
$D(\ce^{(\alpha)})$ formed by  continuous functions with compact support.
A similar identification is not true for the subspaces $B_c(\ce^{(\alpha)})$ of bounded functions unless $\alpha=\infty$.)
Therefore  $(\ce^{(\alpha)})_M=\overline{\ce_0}=\ce^{(\infty)}$ for all $\alpha\in[0,\infty]$.

Secondly, the identification  of  the minimal forms has two distinct cases.
In the limit case $\alpha=\infty$ one has $\ce^{(\infty)}=\overline{\ce_0}$
and $\ce^{(\infty)}$ is regular.
Then  $\ce^{(\infty)}_{m,Y}(\varphi)=\int_Yc\,(\varphi')^2$ for all $\varphi\in D(\overline{\ce_0})$
and each bounded interval $Y$.
Therefore $(\ce^{(\infty)})_m=\overline{\ce_0}=(\ce^{(\infty)})_M$ and the statements of Theorem~\ref{tdrf3.10}
are obviously valid.

Thirdly, if $\alpha \in[0,\infty\rangle$ then $C_c(\ce^{(\alpha)})=C_c(\overline{\ce_0})$ but
$B_c(\ce^{(\alpha)})=B_c(\overline{\ce_0})+\spann(\sigma_+-\sigma_-)
=C_c(\overline{\ce_0})+\spann(\sigma_+-\sigma_-)$.
Therefore the $\ce^{(\alpha)}$ are not inner regular and the conclusions of Theorem~\ref{tdrf3.10} do not apply.
In particular $B_c(\ce^{(\alpha)})\neq B_c(\overline{\ce_0})$.
In fact the forms $(\ce^{(\alpha)})_m$ are not equal to the minimal form $\ce^{(0)}$ of the family $\ce^{(\alpha)}$.
Instead one has  $(\ce^{(\alpha)})_m=\ce^{(\alpha)}$.
   
   \smallskip

    2. If $\nu_+\vee \nu_-\not\in L_\infty(0,1)$ then $H$ has a unique submarkovian extension determined by the Dirichlet form $\ce=\overline{\ce_0}$.
Its domain  is characterized by the continuity conditions $(c\,\varphi')(0_+)=0=(c\,\varphi')(0_-)$.
 In this case $\ce=\ce_M=\ce_m$.
 The semigroup corresponding to the unique Dirichlet form $\ce$ leaves the subspaces $L_2(\Ri_\pm)$ invariant.
 The diffusion  coefficient $c$ is sufficiently degenerate that there is no transmission past the origin which acts as a virtual boundary.
 (Again see \cite{RSi3}, Theorem~1.1.)
 \end{exam}

The foregoing example demonstrates that interior regularity properties are related to the interior degeneracy structure.
Similar effects occur in higher dimensions.
Local degeneracies lead to local discontinuities.

\section{Local extensions}\label{S4}

In this section we examine local Dirichlet forms $\cf$ which extend the local,  inner regular, Dirichlet form $\ce$ or, more generally, extend
 $\ce_M$.
Our aim is to characterize the forms with the property $\ce_m\leq\cf\leq\ce_M$.
Since $\cf$ extends $\ce_M$ it is automatically semi-regular and $\cf\leq \ce_M$. 
The semi-regularity suffices to define the extremal forms $\cf_m,\cf_M$.
Then  $\cf_m\leq\ce_m$ and  $\cf_M\leq\ce_M$.
Thus   the lower bound $\ce_m\leq \cf$ and  its optimality depend on the equality $\ce_m=\cf_m$.
Similarly the optimality of the upper bound $\cf\leq \ce_M$ requires $\cf_M=\ce_M$.
In fact the optimality is related to the inner regularity.

\begin{thm}\label{tdrf4.10}
Let $\ce$ be a local, inner regular,  Dirichlet form and
  $\cf$ a local Dirichlet form  extension   of $\ce_M$.
 
The following conditions are equivalent:
\begin{tabel}
\item\label{tdrf4.10-1}
$\ce_m\leq \cf\leq \ce_M$,
\item\label{tdrf4.10-2}
$B_c(\ce)=B_c(\cf)$.
\end{tabel}
Moreover, if these conditions are satisfied then $\cf$ is inner regular,  $\cf_M=\ce_M$ and $\cf_m=\ce_m$.
\end{thm}
\proof\
\noindent\ref{tdrf4.10-1}$\Rightarrow$\ref{tdrf4.10-2}.$\;$
It follows from the order relation that $B_c(\ce_m)\supseteq B_c(\cf)\supseteq B_c(\ce_M)$.
But  $B_c(\ce_m)=B_c(\ce)=B_c(\ce_M)$ by Theorem~\ref{tdrf3.10}.
Therefore $B_c(\ce)=B_c(\cf)$.

\smallskip

\noindent\ref{tdrf4.10-2}$\Rightarrow$\ref{tdrf4.10-1}.$\;$
First it follows from Condition~\ref{tdrf4.10-2} and Lemma~\ref{inreg} 
that $\cf$ is inner regular.
But $B_c(\cf)=B_c(\ce)$ and $B_c(\ce)=B_c(\ce_M)$ by Theorem~\ref{tdrf3.10}.
Therefore $B_c(\cf)=B_c(\ce_M)$.
Consequently $C_c(\cf)=C_c(\ce_M)$.
Since $C_c(\cf)=C_c(\cf_M)$ it follows that $C_c(\ce_M)=C_c(\cf_M)$.
Therefore $\cf_M=\ce_M$.
Then, however, $\ce_m=(\ce_M)_m=(\cf_M)_m=\cf_m$ by  Theorem~\ref{tdrf3.10} applied first to $\ce$ and then to $\cf$.
The latter application is valid since $\cf$ is inner regular.
One then has $\ce_m=\cf_m\leq \cf\leq \cf_M =\ce_M$.

\smallskip

The last statement of the theorem has been established by the foregoing argument. \hfill$\Box$

\bigskip

If one assumes that the form $\cf$ is  inner regular then it follows by the foregoing argument that 
$\ce_M=\cf_M$ implies $\ce_m=\cf_m$.
Conversely, if  $\ce_m=\cf_m$ then $\ce_M=(\ce_m)_M=(\cf_m)_M=\cf_M$.
But then the conditions $\ce_M=\cf_M$ and $\ce_m=\cf_m$ imply that $\ce_m\leq\cf\leq \ce_M$.
Therefore one has the following conclusion.

\begin{cor}\label{cdrf4.1}
If $\ce$ and $\cf$ are  local, inner regular,  Dirichlet forms with
  $\cf\supseteq \ce_M$  then
the following conditions are equivalent:
\begin{tabel}
\item\label{cdrf4.10-1}
$\ce_m\leq \cf\leq \ce_M$,
\item\label{cdrf4.10-2}
 $\cf_M=\ce_M$,
 \item\label{cdrf4.10-3}
 $\cf_m=\ce_m$.
\end{tabel}
\end{cor}

In fact much more is true under the stronger regularity assumption.

\begin{cor}\label{cdrf4.2}
Assume $\ce$ and $\cf$ are  local, inner regular,  Dirichlet forms with
  $\cf\supseteq \ce_M$.
  Let $S^M$ and $T$ denote the submarkovian semigroups associated with $\ce_M$ and $\cf$, respectively.
Then the following conditions are equivalent:
\begin{tabel}
\item\label{cdrf4.2-1}
$\ce_m\leq \cf\leq \ce_M$,
\item\label{cdrf4.2-2}
$0\leq S^M_t\varphi\leq T_t\varphi$ for all $\varphi\in L_2(X)_+$ and all $t>0$,
\item\label{cdrf4.2-3}
$D(\ce_M)$ is an order ideal of $D(\cf)$,
\item\label{cdrf4.2-4}
$B(\ce_M)$ is an algebraic ideal of $B(\cf)$.
\end{tabel}
\end{cor}
\proof\
The mutual equivalence of the last three conditions follows from Proposition~\ref{tdrf2}
with $\ce$ replaced by $\ce_M$. 
Since $\cf\supseteq \ce_M$  the off-diagonal bounds in Conditions~\ref{tdrf2-2} and \ref{tdrf2-3} of Proposition~\ref{tdrf2}
are automatically fulfilled.
Therefore these conditions reduce to the order ideal property and the algebraic ideal property, respectively.
The equivalences of Proposition~\ref{tdrf2}  are independent of any locality or regularity assumptions.
It  remains to prove equivalence of the last three conditions with the first condition.

\smallskip

\noindent\ref{cdrf4.2-1}$\Rightarrow$\ref{cdrf4.2-2}.$\;$
First Condition~\ref{cdrf4.2-1} is equivalent to the condition $\ce_M=\cf_M$ by Corollary~\ref{cdrf4.1}.
In particular $S^M=T^M$ where $T^M$ denotes the submarkovian semigroup associated with the Dirichlet form $\cf_M$.
But then $0\leq S^M_t\varphi=T^M_t\varphi\leq T_t\varphi$ for all  $\varphi\in L_2(X)_+$ and all $t>0$
by the last statement of Proposition~\ref{pdrf3.11} applied to  $\cf$.
Here it important that $\cf$ is inner regular.

\smallskip
\noindent\ref{cdrf4.2-3}$\Rightarrow$\ref{cdrf4.2-1}.$\;$
Condition~\ref{cdrf4.2-1} is equivalent to the condition $B_c(\ce_M)=B_c(\ce)=B_c(\cf)$ by Theorem~\ref{tdrf4.10}.
Now we argue by negation.

Assume $B_c(\ce_M)\neq B_c(\cf)$.
Since $B_c(\cf)\supseteq B_c(\ce_M)$ it follows that there is a non-zero $\psi\in B_c(\cf)$ such that $\psi\not\in B_c(\ce_M)$.
Replacing $\psi$ by $|\psi|$ and rescaling if necessary we may assume $0\leq \psi\leq 1$.
Since $\ce$ is semi-regular one may then choose $\varphi\in C_c(\ce_M)$ such that $0\leq \varphi\leq 1$ and $\varphi=1$ on $\supp \psi$
by Lemma~\ref{l2}.
Therefore $0\leq \psi\leq \varphi$ and since $\psi\not\in B_c(\ce_M)$ this  contradicts the order ideal property.
\hfill$\Box$

\bigskip

Note that Theorem~\ref{tdrf4.10} incorporates the case of local Dirichlet forms $\cf$ which satisfy $\ce_M\subseteq \cf\subseteq \ce_m$.
This extension/restriction condition immediately implies the ordering condition $\ce_M\geq \cf\geq\ce_m$.
Consequently the conclusions of the theorem are valid.

\smallskip

Example~\ref{ex1d3.1} demonstrates that the conclusions of Theorems~ \ref{tdrf4.10} are essentially optimal.
A similar situation occurs in higher dimensions but then there is a much wider range of possibilities and a greater variety of strongly local extensions.

\begin{exam}\label{exdrf4.2}
Let $\Omega$ be a domain, i.e.\ an open connected set, in $\Ri^d$ and $\ce_0$ the Markovian
form with  $D(\ce_0)=C_c^\infty(\Omega)$ given by
\[
\ce_0(\varphi)=\sum^d_{k,l=1}(\partial_k\varphi, c_{kl}\,\partial_l\varphi)
\]
where $c_{kl}=c_{lk}\in L_{\infty,{\rm loc}}(\Omega)$ and $C=(c_{kl})$ is locally strongly elliptic, i.e.\ for each compact subset $K\subseteq\Omega$ there is a $c_K>0$ such that
$C(x)\geq c_KI$ for (almost) all $x\in\Omega$.
Then $\ce_0$ is closable (see \cite{MR} Section~II.2.b).
Since the truncation $\ce_{0,\xi}$ of $\ce_0$ corresponds to the replacement $c_{kl}\mapsto \xi\,c_{kl}$ the truncated forms with $\xi\geq 0$  are also closable.

Let $\ce_M=\overline \ce_0$.
Then  $\ce_m\,(=(\ce_M)_m)$ is given by 
\begin{equation}
\ce_m(\varphi)=\int_\Omega dx\,\Gamma(\varphi)(x)\;,
\label{ex1}
\end{equation}
where $\Gamma(\varphi)=\sum^d_{k,l=1}c_{kl}\,(\partial_k\varphi)\,(\partial_l\varphi)$, with
 $D(\ce_m)=\{\varphi\in W^{1,2}_{\rm loc}(\Omega): \Gamma(\varphi)+\varphi^2\in L_1(\Omega)\}$.
The forms $\ce_m$ and $\ce_M$ are the extremal forms for each  Dirichlet form extension $\ce$ of $\ce_0$, i.e.\  $\ce_m\leq \ce\leq \ce_M$.
The  form $\ce_M$ is regular, $\ce_m$ is inner regular and   both forms are strongly local.

If $c_{kl}\in W^{1,\infty}_{\loc}(\Omega)$  then $\ce_0$ is the form of the operator $H_0=-\sum^d_{k,l=1}\partial_k\,c_{kl}\,\partial_l$ with $D(H_0)=C_c^\infty(\Omega)$.
Then $H_M$, the operator corresponding to $\ce_M$, is the submarkovian extension of $H_0$ satisfying Dirichlet conditions $\varphi|_{\partial\Omega}=0$  on the boundary $\partial\Omega$.
The operator $H_m$ corresponding to $\ce_m$ formally satisfies Neumann conditions, i.e.\ if $\partial\Omega$ is Lipschitz then $(n.\,C\,\nabla\varphi)|_{\partial\Omega}=0$.
The forms $\ce_F$ defined as   the closure of $\ce_m$  restricted to  $\{\varphi|_\Omega:\varphi\in C_c^\infty(\Ri^d\backslash F)\}$,
where $F$ is a closed subset of $\partial\Omega$, are all inner regular,  strongly local, Dirichlet form extensions of $\ce_0$ which are bounded  above and below by $\ce_M$ and $\ce_m$, respectively.
The  corresponding operators  $H_F$ satisfy a mixture of conditions, Dirichlet on $F$
and Neumann on $\partial\Omega\backslash F$.
The semigroups $S^F$ generated by the $H_F$ dominate the semigroup $S^M$ generated by $H_M$.
In fact since $\ce_F\subseteq \ce_m$ and $D(\ce_F)$ is an order ideal of $D(\ce_m)$ the semigroup $S^m$ generated by $H_m$ dominates the $S^F$.
Thus
$0\leq S^M_t\varphi\leq S^F_t\varphi\leq S^m_t\varphi$ for all $\varphi\geq0$, all $t>0$ and all choices of $F$.

There are also Dirichlet form extensions corresponding to Robin boundary conditions.
If one assumes, for simplicity, that  $\partial\Omega$ is Lipschitz and  sets 
\begin{equation}
\ce^{(\alpha)}(\varphi)=\ce_m(\varphi)+\int_{\partial\Omega}dS\,\alpha\,|\varphi|^2
\label{erob}
\end{equation}
where  $dS$ is the surface measure,  $\alpha$ a positive bounded function on $\partial \Omega$ and $D(\ce^{(\alpha)})$ consists of those $\varphi\in D(\ce_m)$ then the operators 
$H_\alpha$ corresponding to the $\ce^{(\alpha)}$  satisfy the Robin boundary conditions $(n.\,C\,\nabla\varphi+\alpha\,\varphi)|_{\partial\Omega}=0$.
One  has $\ce_m\leq \ce^{(\alpha)}\leq \ce_M$ and $\ce^{(\alpha)}\supseteq \ce_M$ but $\ce_m$ is not an extension of $\ce^{(\alpha)}$.
Nevertheless $\ce_m$ and $\ce_M$ are the extremal forms for the $\ce^{(\alpha)}$.
The construction of $\ce_M$ and $\ce_m$ is by `interior' approximation and the $\alpha$-boundary term plays no role.
The $\ce^{(\alpha)}$ with $\alpha\neq0$ are local, inner regular, forms but they  are not strongly local, because of the boundary term.
One again  has the domination properties $0\leq S^M_t\varphi\leq S^{(\alpha)}_t\varphi\leq S^m_t\varphi$ for all $\varphi\geq0$, all $t>0$ and all choices of $\alpha\geq0$.
The order relation between $S^M$ and $S^{(\alpha)}$ follows from the last statement of Proposition~\ref{pdrf3.11} applied with $S$ replaced by $S^{(\alpha)}$.
The order relation between $S^{(\alpha)}$ and $S^m$ follows from Proposition~\ref{tdrf2} applied to $\ce^{(\alpha)}$ and $\ce_m$.
Since $D(\ce^{(\alpha)})=D(\ce_m)$ the ideal properties of Conditions~\ref{tdrf2-2} and \ref{tdrf2-3} of Proposition~\ref{tdrf2} are evident.
But the off-diagonal bound $\ce^{(\alpha)}(\varphi,\psi)\geq \ce_m(\varphi,\psi)$ for $\varphi,\psi\in D(\ce_m)_+$ follows because $\alpha\geq0$.
These  features also extend to the broader class of Robin conditions  given in \cite{Daners2} \cite{AW2}.
\hfill$\Box$
\end{exam}

This  example illustrates that the locality properties  are quite complex.
If one chooses the strongly local Dirichlet form $\ce_M$ as starting point then the minimal form $\ce_m=(\ce_M)_m$  is strongly local, as are the 
intermediate  forms $\ce_F$, but the Robin forms $\ce^{(\alpha)}$ are only local. 
If, however, one chooses the local Robin form $\ce^{(\alpha)}$ as starting point then the extremal forms are still given by $\ce_m$ and $\ce_M$ and they are both strongly local.
Thus in the latter case there is a strengthening of the  locality property. 

If in Example~\ref{exdrf4.2} the coefficients $c_{kl}\in W^{1,\infty}_{\rm loc}(\Omega)$ then the various self-adjoint operators $H_m$, $H_M$, $ H_F$ and $H^{(\alpha)}$ associated with the forms 
$\ce_m$, $\ce_M$, $\ce_F$ and $\ce^{(\alpha)}$ are all submarkovian extensions of the symmetric operator $H_0$.
In particular $H_m$ is the smallest such extension and $H_M$ the largest, e.g.\ $H_m\leq H_F\leq H_M$ and $H_m\leq H^{(\alpha)}\leq H_M$.
Therefore $H_0$ is Markov unique, i.e.\ it has a unique submarkovian extension, if and only if $H_m=H_M$ or, equivalently, $\ce_m=\ce_M$.
In the next section we examine this latter uniqueness criterion in the general situation.

\section{Applications}\label{S5}

In this section we consider several implications of the foregoing results.
We concentrate on two topics,
the dependence  on $\ce$ of the set-theoretic distance function  
$d^{\,(\ce)}(\,\cdot\,;\,\cdot\,)$ canonically associated with the form
and the uniqueness criterion $\ce_m=\ce_M$.
In contrast to the earlier discussion strong locality is now essential for  much of the analysis.
In fact 
strong locality of the form $\ce_m$ is of prime importance.
This implies strong locality of $\ce_M$, because $\ce_m\supseteq \ce_M$,
but it does not necessarily require strong locality of $\ce$
(see Remark~\ref{rdrf3.1}  and the
discussion following Example~\ref{exdrf4.2}).
The second key ingredient in the discussion  is the algebraic ideal property, i.e.\  the observation that $B(\ce_M)$ is an algebraic
ideal of $B(\ce_m)$.

\subsection{Distances}\label{S5.2}

One can associate with  a general strongly local Dirichlet form  $\ce$  a set-theoretic distance function \cite{Stu2} \cite{HiR} \cite{AH} \cite{ERSZ2}.
In particular the distance as  defined by Ariyoshi and Hino \cite{AH}  for a strongly local Dirichlet form $\ce$ on the  measure space $X,\mu$ is independent
of any topology. 
It is a positive function $d^{\,(\ce)}(A\,;B)$ over pairs of measurable sets $A, B$ with $\mu(A),\mu(B)\in\langle0,\infty\rangle$.
We next sketch the construction of this function and refer to \cite{AH} for details.

First for each  closed subset $A$  of $X$ set
\[
D_{\!A}(\ce)=\{\varphi\in D(\ce): \supp\varphi= A\}
\]
and $B_{\!A}(\ce)=D_A(\ce)\cap L_\infty(X)$.
Secondly, following  \cite{AH} Definition~2.1, define a nest as an increasing family  $\ca=\{A_\lambda\}_{\lambda>0}$ of sets $A_\lambda$ with finite $\ce$-capacity such that $D_{\!\ca}(\ce)=\bigcup_{\lambda>0}D_{\!A_\lambda}(\ce)$ is a core of $\ce$ or, equivalently,  that $B_{\!\ca}(\ce)=\bigcup_{\lambda>0}B_{\!A_\lambda}(\ce)$ is a core.
It is not  evident that nests of this type exist but this follows from  \cite{AH} Lemma~3.1 by an argument based on the theory of excessive functions.

Thirdly  define the local subspace of $D(\ce)$ corresponding to the nest $\ca$ by
\[
D_{\!\ca; \rm loc}(\ce)= \{\varphi\in M(X): \mbox{ there exist } \varphi_\lambda\in D(\ce) \mbox{ such that } \varphi=\varphi_\lambda \mbox{ on } A_\lambda\}
\]
and $B_{\!\ca; {\rm loc}}(\ce)=D_{\!\ca; \rm loc}\cap L_\infty(X)$ (see \cite{AH} Definition~2.4).
Since the  $A_\lambda$ are sets of   finite $\ce$-capacity it immediately follows that $\one \in B_{\!\ca; {\rm loc}}(\ce)$.
Moreover, if $F$ is a normal contraction and $\varphi\in D_{\!\ca; \rm loc}(\ce)$ then $F\circ\varphi\in D_{\!\ca; \rm loc}(\ce)$ since $\ce$ is a Dirichlet form.

Fourthly if  $\ce$  is 
local one can extend the earlier definition of the truncated forms to the local functions.
It is convenient for the sequel to adopt the notation of  \cite{ERSZ2} and set 
\[
\ci^{(\ce)}_{\;\varphi}(\xi)=\ce_\xi(\varphi)=\ce(\varphi,\xi\,\varphi)-2^{-1}\ce(\xi,\varphi^2)
\]
for $\xi\in B(\ce)_+$   and $\varphi\in D(\ce_\xi)=B(\ce)$.
It then follows from the 
locality of $\ce$ that $\ci^{(\ce)}_{\;\varphi}(\xi)=0$ if $\varphi\,\xi=0$.
Consequently  for each $\varphi\in B_{\!\ca; {\rm loc}}(\ce)$ one can introduce  the  form $\hat\ci^{(\ce)}_{\;\varphi}$  on $B_{\!\ca}(\ce)_+$
by setting
\[
\hat\ci^{(\ce)}_{\;\varphi}(\xi)=\ci^{(\ce)}_{\;\varphi_\lambda}(\xi)
\]
for $\xi\in B_{\!A_\lambda}(\ce)_+$ 
and for any $\varphi_\lambda\in D(\ce)$ with  $\varphi|_{A_\lambda}=\varphi_\lambda$.

Finally  define
\[
|||\,\hat\ci^{(\ce)}_{\;\varphi}|||=\sup\{\,\hat\ci^{(\ce)}_{\;\varphi}(\xi) : \xi\in B_{\!\ca}(\ce)_+\,,\, \|\xi\|_1\leq1\}
\]
and  introduce
\[
D_{0,\ca}(\ce)=\{\varphi\in B_{\!\ca; {\rm loc}}(\ce): |||\,\hat\ci^{(\ce)}_{\;\varphi}|||\leq1\,\}
\;.
\]
This corresponds to Definition~2.6 of \cite{AH}.
It then  follows from  Proposition~3.9 of \cite{AH} that  if $\ce$ is strongly local then $D_{0,\ca}(\ce)$  is independent of the particular choice of nest, i.e.\ $D_{0,\ca}(\ce)=D_{0,\cb}(\ce)$
for any pair of nests $\ca,\cb$.
Hence one may  simplify notation by setting $D_0(\ce)=D_{0,\ca}(\ce)$.
Then for each pair of  measurable sets $A,B\subset X$ with $\mu(A), \mu(B)>0$ 
 the distance between $A$ and $B$ corresponding to $\ce$ is defined by
\[
d^{\,(\ce)}(A\,;B)=\sup\Big\{\essinf_{x\in A}\varphi(x)-\esssup_{y\in B}\varphi(y): \varphi\in D_0(\ce)\Big\}
\;.
\]
This definition agrees with that of Ariyoshi and Hino and we 
use their results in the sequel.

At this point we return to the setting of Sections~\ref{S3} and \ref{S4} and compare the distances corresponding to 
extensions of a strongly local inner regular form $\ce$.
The algebraic ideal properties are of importance.
Recall that if $\ce$ and $\cf$ are two Dirichlet forms with $\cf\geq \ce$ then $D(\cf)\subseteq D(\ce)$.
Therefore $B_{\!\ca}(\cf)\subseteq B_{\!\ca}(\ce)$ and $B_{\!\ca; {\rm loc}}(\cf)\subseteq B_{\!\ca; {\rm loc}}(\ce)$.
Further if $B(\cf)$ is  an algebraic ideal of $B(\ce)$, i.e.\  if $B(\cf)\,B(\ce)\subseteq B(\cf)$, then
 $B_{\!\ca; {\rm loc}}(\cf)\,B_{\!\ca}(\ce)\subseteq B_{\!\ca; {\rm loc}}(\cf)$.

The principal result concerns the extremal forms. 

\begin{thm}\label{tcap1} Let $\ce$ be an inner regular, strongly local, Dirichlet form and $\ce_m$, $\ce_M$ 
the corresponding minimal and maximal Dirichlet forms.
Then
\[
d^{\,(\ce)}(A\,;B)=d^{\,(\ce_M\!)}(A\,;B)
\]
for all measurable subsets $A,B$ with $\mu(A),\mu(B)\in\langle0,\infty\rangle$.
Moreover, if $\ce_m$ is strongly local then one also has
\[
d^{\,(\ce_m\!)}(A\,;B)=d^{\,(\ce)}(A\,;B)
\]
for all  $A,B$ with $\mu(A),\mu(B)\in\langle0,\infty\rangle$.
\end{thm}
\proof\
First  strong locality of $\ce_M$ follows from strong locality of $\ce$ since $\ce_M\subseteq \ce$.
Therefore $d^{\,(\ce)}$ and $d^{\,(\ce_M\!)}$ are both well-defined.
Now we begin by proving  that $d^{\,(\ce_M\!)}(A\,;B)\leq d^{\,(\ce)}(A\,;B)$ for all $A,B$.
The proof is a variation of the argument used to establish Proposition~5.3 in \cite{ERSZ2}.

\smallskip

Fix $\varphi\in B_{\!\ca;{\rm loc}}(\ce_M)_+$ with  $|||\,\hat\ci^{(\ce_M\!)}_{\;\varphi}|||<\infty$.
Now $\one\in  B_{\!\ca;{\rm loc}}(\ce_M)_+$ and if  $\xi\in B_{\!A_\lambda}(\ce)_+$ then
\begin{eqnarray*}
\hat\ci^{(\ce)}_{\;\one+\varphi}(\xi)&=&\ce((\varphi_\lambda+\one_\lambda),\xi\,(\varphi_\lambda+\one_\lambda))-2^{-1}\ce(\xi,(\varphi_\lambda+\one_\lambda)^2)\\[5pt]
&=&\ce(\varphi_\lambda,\xi\,\varphi_\lambda) +\ce(\varphi_\lambda,\xi)-2^{-1}\,(\ce(\xi,\varphi_\lambda^2)+2\,\ce(\xi,\varphi_\lambda))=
\hat\ci^{(\ce)}_{\;\varphi}(\xi)
\end{eqnarray*}
by strong locality of $\ce$.
Therefore  $\hat\ci^{(\ce)}_{\;\varphi}(\xi)=\hat\ci^{(\ce)}_{\;\one+\varphi}(\xi)$ for all $\xi\in B_{\!\ca}(\ce)_+$.

Next $(\one+\varphi)^{1/2}=\one+F\circ\varphi$ with $F(x)=(1+|x|)^{1/2}-1$.
But $F$ is a normal contraction.
Therefore $(\one+\varphi)^{1/2}\in B_{\!\ca;{\rm loc}}(\ce_M)_+$.
Now the polarized form
\[
\ci^{(\ce)}_{\;\varphi,\psi}(\xi)=\ce_{\xi}(\varphi,\psi)=2^{-1}\Big(\ce(\varphi, \xi\,\psi)+\ce(\xi\,\varphi, \psi)-\ce(\xi,\varphi\,\psi)\Big)
\]
satisfies  the Leibniz relation
\[
\ci^{(\ce)}_{\;\varphi_1\varphi_2,\psi}(\xi)=\ci^{(\ce)}_{\;\varphi_1,\psi}(\varphi_2\,\xi)+\ci^{(\ce)}_{\;\varphi_2,\psi}(\varphi_1\,\xi)
\]
as a consequence of strong locality.
It then readily follows that
\begin{equation}
\hat\ci^{(\ce)}_{\;\varphi}(\xi)=\hat\ci^{(\ce)}_{\;\one+\varphi}(\xi)=4\,\hat\ci^{(\ce)}_{\;(\one+\varphi)^{1/2}}((\one+\varphi)\,\xi)
\label{eess}
\end{equation}
for all $\xi\in B_{\!\ca}(\ce)_+$.
But $\eta=(\one+\varphi)\,\xi\in B_{\!\ca;{\rm loc}}(\ce_M)_+$ by the ideal property.
Moreover,
$\hat\ci^{(\ce)}_{\;(\one+\varphi)^{1/2}}(\eta)\leq \hat\ci^{(\ce_M\!)}_{\;(\one+\varphi)^{1/2}}(\eta)$ for all 
$\eta\in B_{\!\ca}(\ce_M)_+$ by the order property  $\ce\leq \ce_M$.
Therefore 
\[
\hat\ci^{(\ce)}_{\;\varphi}(\xi)\leq 4\,\hat\ci^{(\ce_M\!)}_{\;(\one+\varphi)^{1/2}}((\one+\varphi)\,\xi)
\]
for all $\xi\in B_{\!\ca}(\ce)_+$.
It immediately follows that 
\[
|||\,\hat\ci^{(\ce)}_{\;\varphi}|||\leq 4\,(1+\|\varphi\|_\infty)\,|||\,\hat\ci^{(\ce_M)}_{\;(\one+\varphi)^{1/2}}|||
\;.
\]
Then, however, one deduces from (\ref{eess})  with $\ce$ replaced by $\ce_M$ that
\[
4\,\hat\ci^{(\ce_M\!)}_{\;(\one+\varphi)^{1/2}}((\one+\varphi)\,\eta)=\hat\ci^{(\ce_M\!)}_{\;\varphi}(\eta)
\]
for all $\eta\in B_{\!\ca}(\ce_M)_+$.
Therefore 
\[
|||\,\hat\ci^{(\ce_M\!)}_{\;(\one+\varphi)^{1/2}}|||\leq 4^{-1}\,|||\,\hat\ci^{(\ce_M\!)}_{\;\varphi}|||
\;.
\]
Combining these estimates one concludes that 
\[
|||\,\hat\ci^{(\ce)}_{\;\varphi}|||\leq (1+\|\varphi\|_\infty)\,|||\,\hat\ci^{(\ce_M\!)}_{\;\varphi}|||
\]
for all $\varphi\in B_{\!\ca;{\rm loc}}(\ce_M)_+$.

Finally, replacing $\varphi$ by $\tau\,\varphi$ with $\tau>0$ and  noting that $|||\,\hat\ci^{(\ce)}_{\;\tau \varphi}|||=\tau^2\,|||\,\hat\ci^{(\ce)}_{\;\varphi}|||$
and $|||\,\hat\ci^{(\ce_M\!)}_{\;\tau \varphi}|||=\tau^2\,|||\,\hat\ci^{(\ce_M\!)}_{\;\varphi}|||$ one deduces that 
\[
|||\,\hat\ci^{(\ce)}_{\;\varphi}|||\leq (1+\tau\,\|\varphi\|_\infty)\,|||\,\hat\ci^{(\ce_M\!)}_{\;\varphi}|||
\]
for all $\varphi\in B_{\!\ca;{\rm loc}}(\ce_M)_+$.
Therefore in the limit $\tau\to0$ one has $|||\,\hat\ci^{(\ce)}_{\;\varphi}|||\leq|||\,\hat\ci^{(\ce_M\!)}_{\;\varphi}|||$.
Consequently   $d^{\,(\ce_M\!)}(A\,;B)\leq d^{\,(\ce)}(A\,;B)$ for all $A,B$.

The second step in the proof is  to establish  the converse inequalities.
But this is a corollary of  the small time asymptotics of the  semigroups  $S$  and $S^M$
established by Ariyoshi and Hino \cite{AH} and the order
property of $S$  and $S^M$ given  by Proposition~\ref{pdrf3.11}.

First it follows from \cite{AH},  Theorem~2.7, that 
\[
d^{\,(\ce)}(A\,;B)^2=-\lim_{t\to0}4\,t\log(\one_A,S_t\one_B)\;\;\;\;\;{\rm and}\;\;\;\;\;
d^{\,(\ce_M\!)}(A\,;B)^2=-\lim_{t\to0}4\,t\log(\one_A,S^M_t\one_B)
\]
for for all $A,B$ with $\mu(A),\mu(B)\in\langle0,\infty\rangle$.
(Note that our convention for the semigroup generator differs from that of \cite{AH} by a factor $2$.)

Secondly, it follows from Proposition~\ref{pdrf3.11} that 
\[
(\one_A,S^M_t\one_B)\leq (\one_A,S_t\one_B)
\]
for all $t>0$.
Therefore 
\[
d^{\,(\ce_M\!)}(A\,;B)^2=-\lim_{t\to0}4\,t\log(\one_A,S^M_t\one_B)
\geq -\lim_{t\to0}4\,t\log(\one_A,S_t\one_B)=d^{\,(\ce)}(A\,;B)^2
\;.
\]
and one concludes that $d^{\,(\ce_M\!)}(A\,;B)\geq d^{\,(\ce)}(A\,;B)$ for all $A,B$ with $\mu(A),\mu(B)\in\langle0,\infty\rangle$.

\smallskip

The proof of the second statement of the theorem, i.e.\ the equality  $d^{\,(\ce_m\!)}(A\,;B)=d^{\,(\ce)}(A\,;B)$, follows by repeating the
foregoing argument with $\ce$ replaced by $\ce_m$ to conclude that $d^{\,(\ce_m\!)}(A\,;B)=d^{\,(\ce_M\!)}(A\,;B)$ and then combining
this with the first statement of the theorem.
In this case the   Ariyoshi--Hino asymptotic estimate for $S^m$  requires  the  strong locality of $\ce_m$.
Moreover, it uses the ordering for $S^m$ and $S^M$  given by Proposition~\ref{pdrf3.14}
in place of the ordering of $S$ and $S^M$ given by Proposition~\ref{pdrf3.11}.
\hfill$\Box$

\begin{cor}\label{ccap1}
Let  $\ce$ and $\cf$ be strongly local inner regular forms with $\ce_M\subseteq\cf$.
Assume  $\ce_m\leq \cf$.
Then
\[
d^{\,(\cf)}(A\,;B)=d^{\,(\ce)}(A\,;B)=d^{\,(\ce_M\!)}(A\,;B)=d^{\,(\cf_M)}(A\,;B)
\]
for all measurable subsets $A,B$ with $\mu(A),\mu(B)\in\langle0,\infty\rangle$.
Moreover, if $\ce_m$ is strongly local one also has $d^{\,(\cf_m\!)}(A\,;B)=d^{\,(\ce)}(A\,;B)$.
\end{cor}
\proof\
It follows from Theorem~\ref{tcap1} applied to $\cf$ that 
\[
d^{\,(\cf)}(A\,;B)=d^{\,(\cf_M\!)}(A\,;B)
\]
for all $A,B$ with $\mu(A),\mu(B)\in\langle0,\infty\rangle$.
But  $\ce_m\leq \cf\leq \ce_M$ by assumption and Corollary~\ref{cdrf4.1}
 establishes that this order relation  is equivalent to the identity $\cf_M=\ce_M$. 
 Therefore 
 \[
d^{\,(\cf_M\!)}(A\,;B)=d^{\,(\ce_M\!)}(A\,;B)=d^{\,(\ce)}(A\,;B)
\]
 with the second identity following from
Theorem~\ref{tcap1}.
The last statement of the corollary follows by similar reasoning.
\hfill$\Box$

\bigskip

Note that although the distance $d^{\,(\ce)}(A;B)$ can in principle be infinite this is not the case if $S$ is irreducible, i.e.\ if $(\one_A, S_t\one_B)>0$
for all $A, B$ with $\mu(A),\mu(B)\in \langle0,\infty\rangle$ and $t>0$.
This follows because $S$ satisfies the Davies--Gaffney bounds $(\one_A, S_t\one_B)\leq \mu(A)^{1/2}\,\mu(B)^{1/2}\,e^{-d^{\,(\ce)}(A;B)^2/4t}$ for all $t>0$ (see \cite{AH}
Theorem~4.1).
On the other hand it is possible that $d^{\,(\ce)}$ is identically zero.
In particular this happens if $D_0(\ce)=\{0\}$, or  $D_0(\ce)=\{\lambda\one_X\}$ if $\mu(X)<\infty$.
This follows because $\varphi\in D_0(\ce)$ requires the corresponding Radon measure $\mu_\varphi$ to be absolutely continuous
with respect to $\mu$ and  $d\mu_\varphi/d\mu\leq 1$.
This condition typically fails for diffusions on fractals equipped with the Hausdorff measure.

\smallskip

Next we illustrate the foregoing results with the example of strictly elliptic forms.
The equality of the various distances is a statement of independence of the geometry
from the choice of boundary conditions for the corresponding diffusion process.

\begin{exam}\label{exdrf5.1}
Let $\ce_0$ denote the Markovian form of Example~\ref{exdrf4.2} on $L_2(\Omega)$ where $\Omega$ is a domain in $\Ri^d$.
Let $\ce_M$  and  $\ce_m (=(\ce_M)_m)$ denote the corresponding extremal Dirichlet forms.
Then  $\ce_m$,  $\ce_M$ and all Dirichlet forms $\ce$  with $\ce_m\supseteq \ce\supseteq \ce_M$ are  strongly local.
Therefore $d^{\,(\ce_m\!)}(A\,;B)=d^{(\ce)}(A\,;B)=d^{\,(\ce_M\!)}(A\,;B)$ for all measurable subsets $A,B$ with $0<|A|,|B|<\infty$.
Now  we argue that if  the coefficients $c_{kl}$ are Lipschitz continuous then $d^{\,(\ce_M\!)}(A\,;B)$ is the geodesic distance
between the  open  subsets $A$ and $B$ corresponding to the metric $C^{-1}$.

First  since $\ce_M$ is regular one can compute the set-theoretic distance with the nest of  compact subsets of $\Omega$.
But
\[
\ci^{(\ce_M\!)}_{\;\psi}(\xi)=\ce_M(\psi,\xi\,\psi)-2^{-1}\ce_M(\xi,\psi^2)=\int_\Omega dx\,\xi(x)\,\Gamma(\psi)(x)
\]
for $\xi, \psi \in D(\ce_M)$  with compact support and with  $\xi\geq0$.
Now if  $\psi\in B_{{\rm loc}}(\ce_M)$ and $K\subset \Omega$ is  compact
one can choose $\hat\psi\in B(\ce_M)$ such that $\psi|_K=\hat\psi|_K$.
Then  $\hat\ci^{(\ce_M\!)}_{\;\psi}(\xi)=\ci^{(\ce_M\!)}_{\;\hat\psi}(\xi)$ for all $\xi$ with $\supp\xi\subseteq K$.
Therefore if $|||\,\hat\ci^{(\ce_M\!)}_{\;\psi}|||\leq1$  one has 
\[
\Big|\int_K dx\,\xi(x)\,\Gamma(\psi)(x)\Big|=\Big|\int_Kdx\,\xi(x)\,\Gamma(\hat\psi)(x)\Big|=|\,\ci^{(\ce_M\!)}_{\;\hat\psi}(\xi)|
=|\,\hat\ci^{(\ce_M\!)}_{\psi}(\xi)|\leq \|\xi\|_1
\;.
\]
Hence $\sup_{x\in K}|\,\Gamma(\psi)(x)|\leq 1$ uniformly for all $K$ and all $\psi\in B_{{\rm loc}}(\ce_M)$.
Thus $D_0(\ce_M)=\{\psi\in W^{1,\infty}_{\rm loc}(\Omega):\|\,\Gamma(\psi)\|_\infty\leq 1\}$.
Therefore $d^{\,(\ce_M\!)}(A\,;B)=\inf_{x\in A,y\in B}d_C(x\,;y)$ where 
\[
d_C(x\,;y)
=\sup\{\psi(x)-\psi(y): \psi\in W^{1,\infty}_{\rm loc}(\Omega)\,,\, \|\Gamma(\psi)\|_\infty\leq 1\}
\;.
\]
But the latter expression is one of the well known characterizations of the geodesic distance.

The forms $\ce_m, \ce, \ce_M$ are distinguished by different boundary conditions and the foregoing calculation establishes that the 
set-theoretic distance is independent of the boundary conditions.
This conclusion also follows for the strongly local forms $\ce_F$ corresponding to Dirichlet boundary conditions on the closed subset $F$
of the boundary and Neumann conditions on the complement $\partial\Omega\backslash F$ since $\ce_m\leq \ce_F\leq \ce_M$.

Finally  if $\Omega$ has a smooth boundary one can also define the  forms $\ce^{(\alpha)}$ with Robin boundary conditions by (\ref{erob}).
These forms are local but not strongly local.
Therefore the Ariyoshi--Hino definition of the set-theoretic distance is not applicable.
Nevertheless one can deduce that the Robin semigroups have the same small time asymptotic behaviour as the Dirichlet and Neumann semigroups.
This follows because
\[
(\one_A, S^M_t\one_B)\leq (\one_A, S^{(\alpha)}_t\one_B)\leq (\one_A, S^m_t\one_B)
\]
for all $A, B\subseteq \Omega$ with $0<|A|,|B|<\infty$ and $t>0$ by the discussion in  Example~\ref{exdrf4.2}. 
Then since $d^{\,(\ce_m\!)}(A\,;B)=d^{\,(\ce_M\!)}(A\,;B) \,(=d^{(\ce)}(A\,;B))$  it follows by the Ariyoshi--Hino asymptotic estimates
for $S^M$ and $S^m$ that
\[
d^{(\ce)}(A\,;B)^2=-\lim_{t\to0}4\,t\log (\one_A, S^{(\alpha)}_t\one_B)
\;.
\]
A similar argument is valid for the Robin semigroups defined for arbitrary domains in 
\cite{AW2} sinnce the sandwich estimate for the semigroups is established in \cite{AW}.
\end{exam}

\subsection{Uniqueness}\label{S5.1}

In this subsection we consider the condition $\ce_m=\ce_M$.
This condition has been extensively analysed as a criterion of Markov 
uniqueness of second-order elliptic operators on domains $\Omega$ in $\Ri^d$.
In this latter setting it is equivalent to several  rather different conditions
(see \cite{Ebe}, \cite{RSi4} and \cite{Rob8} for further details and references).
In particular it known to be equivalent to two distinct types of capacity condition.
The first result of this nature was due to Maz'ya (see \cite{Maz} Section~2.7 or \cite{Ebe} Theorem~3.6).
Our next aim is to demonstrate that a  characterization similar to that of Maz'ya  can be established in the 
general Dirichlet form setting.

Recall that  $B_{\capp}(\ce)$ denotes the subspace of bounded functions in $D(\ce)$ whose supports have finite $\ce$-capacity.
 If $A$ is a subset of $X$ with finite $\ce$-capacity we then set $D_{\capp, A}(\ce)=\{\psi\in D_{\capp}(\ce):\supp\psi\subseteq A\}$
 and $B_{\capp, A}(\ce)=D_{\capp, A}(\ce)\cap L_\infty(X)$.
 
\begin{prop}\label{pcap3.1} Let $\ce$ be a  local, inner regular,  Dirichlet form and $\ce_m$, $\ce_M$ the corresponding extremal forms.
Consider the  following conditions:
\begin{tabelpairs}
 \item\label{pcap3.1-3}
 for each subset $A$  of $X$  with  finite $\ce_m$-capacity
 there exists a sequence  $\{\eta_n\}_{n\geq1}$ of \,$\eta_n \in C_c(\ce_M)$ with $0\leq\eta_n\leq 1$ $($of\, $\eta_n \in B(\ce_M))$   such that
\[
\lim_{n\to\infty}\|(\one_X-\eta_n)\,\varphi\|_{D(\ce_m)}=0
\]
 for all  $\varphi \in B_{\capp, A}(\ce_m)$,
\end{tabelpairs}
\begin{tabel}
 \setcounter{teller}{1}
 \item\label{pcap3.1-2}
$\hspace{8mm}\ce_m=\ce_M$.
\end{tabel}
  
 Then {\rm \ref{pcap3.1-3}$\Rightarrow$\ref{pcap3.1-3}$^\prime\!\!\Rightarrow$\ref{pcap3.1-2}}.
Moreover, if $\ce_M$ is strongly local then {\rm \ref{pcap3.1-2}$\Rightarrow$\ref{pcap3.1-3}} and the three conditions are equivalent.
 \end{prop}
 \proof\
\noindent\ref{pcap3.1-3}$\Rightarrow$\ref{pcap3.1-3}$^\prime$.$\;$
This is evident.

\smallskip

\noindent\ref{pcap3.1-3}$^\prime\!\!\Rightarrow$\ref{pcap3.1-2}.$\;$
Fix $\varphi \in B_{\capp}(\ce_m)$ and set $A=\supp\varphi$.
Let $\eta_n$ be the sequence in Condition~\ref{pcap3.1-3}$^\prime$ corresponding to $A$ and set $\varphi_n=\eta_n\, \varphi$.
Since  $B(\ce_M)$ is an algebraic ideal  of $B(\ce_m)$, by Proposition~\ref{pdrf3.14},  it follows that $\varphi_n\in B(\ce_M)$.
But $\ce_m\supseteq \ce_M$ by Theorem~\ref{tdrf3.10}.\ref{tdrf3.10-3}.
Therefore Condition~\ref{pcap3.1-3}$^\prime$ implies that the sequence $\varphi_n$ is convergent with respect to the $D(\ce_M)$-graph norm. 
This establishes that  $\varphi\in B(\ce_M)$.
Hence  $B_{\capp}(\ce_m)\subseteq  B(\ce_M)$ and $\ce_M(\varphi)=\ce_m(\varphi)$
for all $\varphi\in B_{\capp}(\ce_m)$.
But $B_{\capp}(\ce_m)$ is a core of $\ce_m$ by  Proposition~\ref{pcap1}.
Therefore $\ce_m=\ce_M$.

\smallskip

\noindent\ref{pcap3.1-2}$\Rightarrow$\ref{pcap3.1-3}.$\;$
We now assume $\ce_M$. 
Hence $\ce_m$ is strongly local by Condition~\ref{pcap3.1-2}.
Since $A$ has  finite $\ce_m$-capacity  there is an $\eta\in D(\ce_m)$
with $\eta=1$ on $A$.
Then $\eta\,\varphi=\varphi$ for all $\varphi \in B_{\capp, A}(\ce_m)$.
But  $\eta\in D(\ce_M)$,  by Condition~\ref{pcap3.1-2}, and it
 follows from the Dirichlet property that one may assume that $0\leq \eta\leq 1$.
Therefore one can choose a sequence $\eta_n\in C_c(\ce_M)$ with $0\leq\eta_n\leq1$ such that $\|\eta_n-\eta\|_{D(\ce_M)}\to0$ as $n\to\infty$.
In particular $\sup_n\ce_M(\eta_n)<\infty$.

Next it follows that $\|(\eta_n-\eta)\varphi\|_2\leq \|\eta_n-\eta\|_2\|\varphi\|_\infty\to0$ as $n\to \infty$
for all $\varphi \in    B_{\capp, A}(\ce_m)$.
Moreover,
\begin{eqnarray*}
\ce_m((\eta_n-\eta)\varphi)
\leq 2\,\ce_m(\varphi)+2\,\ce_m(\eta_n\,\varphi)
\;.
\end{eqnarray*}
But it follows from the strong locality of $\ce_m$, by a straightforward application of Theorem~5.2.1 in \cite{BH} that 
\[
\ce_m(\eta_n\,\varphi)\leq 2\, (\ce_m)_{\eta_n^2}(\varphi)+  2\, (\ce_m)_{\varphi^2}(\eta_n)
\;.
\]
Therefore one deduces that 
\[
\ce_m((\eta_n-\eta)\varphi)\leq  6\,\ce_m(\varphi)+4\,\ce_M(\eta_n)\|\varphi\|_\infty^2
\;.
\]
Consequently $\sup_n\ce_m((\eta_n-\eta)\varphi)<\infty$.
Hence there is a subsequence $(\eta_{n_k}-\eta)\,\varphi$ which is weakly convergent  to zero in the Hilbert space
$D(\ce_m)$ equipped with the graph norm and  with Cesaro mean strongly convergent to zero.
Then replacing $\eta_n$ by $\eta'_n=n^{-1}\sum_{k=1}^n\eta_{n_k}$ one deduces that
$\|(\one_X-\eta_n')\varphi\|_{D(\ce_m)}=\|(\eta-\eta'_n)\,\varphi\|_{D(\ce_m)}\to0$.
Therefore the sequence $\eta'_n$ satisfies Condition~\ref{pcap3.1-3}. 
\hfill$\Box$

\bigskip

The situation is simpler if $\mu(X)<\infty$ and $\ce$ is strongly local.
Then the identity function $\one_X$ is in the domain of the form $\ce_m$ and $\ce_m(\one_X)=0$.
This follows by first noting that  $\one_X\in L_2(X)$.
Secondly, if  $\xi, \varphi\in C_c(\ce)$ with $0\leq \xi\leq 1$ and $\varphi=1$ on $\supp\xi$ then
\[
\ce_\xi(\varphi)=\ce(\varphi,\xi\varphi)-2^{-1}\ce(\xi,\varphi^2)=\ce(\varphi,\xi)-2^{-1}\ce(\xi,\varphi^2)=0
\]
by strong locality of $\ce$.
Therefore  $\one_X\in D(\ce_{m,Y;0})$ and $\ce_{m,Y;0}(\one_X)=0$ for all bounded open subsets $Y$ of $X$.
Then by the definition of $\ce_m$ one has $\one_X\in D(\ce_m)$ and $\ce_m(\one_X)=0$.

\begin{cor}\label{ccap3.1} Let $\ce$ be a strongly local, inner regular,  Dirichlet form and $\ce_m$, $\ce_M$ the corresponding extremal forms.
Assume $\mu(X)<\infty$.
Then the  following conditions are equivalent:
\begin{tabelpairs}
\item\label{ccap3.1-2}
  there exists a sequence
 $\{\eta_n\}_{n\geq1}$ of \,$\eta_n\in C_c(\ce_M)$  with $0\leq \eta_n\leq 1$  $($of \,$\eta_n\in  B(\ce_M))$  such that 
$\lim_{n\to\infty}\|\one_X-\eta_n\|_{D(\ce_m)}=0$,
\end{tabelpairs}
\begin{tabel}
 \setcounter{teller}{1}
\item\label{ccap3.1-1}
$\hspace{8mm}\ce_m=\ce_M$.
\end{tabel}
\end{cor}
\proof\   \ref{ccap3.1-2}$\Rightarrow$\ref{ccap3.1-2}$^\prime$.$\;$ This is evident.

\smallskip

\noindent \ref{ccap3.1-2}$^\prime$$\Rightarrow$\ref{ccap3.1-2}.$\;$ One may first choose $\xi_n\in B_c(\ce_M)$ such that $\|\eta_n-\xi_n\|_{D(\ce_M)}\leq n^{-1}$ for all $n\geq1$.
 Secondly one may choose $\zeta_n\in C_c(\ce_M)$ such that $\|\xi_n-\zeta_n\|_{D(\ce_M)}\leq n^{-1}$ for all $n\geq1$ by inner regularity.
 Then one has $\|\one_X-\zeta_n\|_{D(\ce_m)}\leq \|\one_X-\eta_n\|_{D(\ce_m)}+2\,n^{-1}$.
  Finally it follows from the Dirichlet property that the sequence  $0\vee\zeta_n\wedge 1$ satisfies Condition~\ref{ccap3.1-2}.
 
 \smallskip
 
  \noindent  \ref{ccap3.1-2}$\Rightarrow$\ref{ccap3.1-1}.$\;$ 
 Condition~\ref{ccap3.1-2} implies that $\one_X\in B(\ce_M)$. 
But $B(\ce_M)$ is an algebraic ideal of $B(\ce_m)$ by Proposition~\ref{pdrf3.14}.
Therefore  $\varphi=\one_X\varphi\in B(\ce_M)$ for all $\varphi\in B(\ce_m)$.
Consequently $B(\ce_m)=B(\ce_M)$ and  $\ce_m=\ce_M$.

\smallskip

\noindent\ref{ccap3.1-1}$\Rightarrow$\ref{ccap3.1-2}.$\;$  If $\ce_m=\ce_M$ then $C_c(\ce_M)$ is a core of $\ce_m$.
Moreover,  $\ce_M=\ce$ is strongly local.
 But $\one_X\in D(\ce_m)$ by the preceding discussion.
 Therefore Condition~\ref{ccap3.1-2} follows immediately from Proposition~\ref{pcap3.1}.
\hfill$\Box$

\bigskip

Mazya's criterion for uniqueness  \cite{Maz}, Section~2.7, was originally formulated for second-order elliptic operators
in divergence form on $\Ri^d$.
Then the capacity estimates of Proposition~\ref{pcap3.1} and Corollary~\ref{ccap3.1} give bounds on the possible growth
of the coefficients at infinity which ensure that the `boundary at infinity'  is inaccessible to the corresponding diffusion.
Complementary estimates have been given for operators on domains with boundaries in \cite{RSi4} \cite{RSi5}.
We conclude by establishing these latter estimates for bounded domains.

First note that the $\ce$-capacity of a subset $A\subseteq X$ is defined by
\[
\capp_\ce(A)=\inf\{\|\varphi\|_{D(\ce)}^2:\varphi\in D(\ce), \, \varphi=1 \mbox{  on } A\,\}
\]
and one can restrict the infimum to $\varphi$ with $0\leq\varphi\leq1$ by the Dirichlet property of $\ce$.
Secondly let $\Omega$ be a domain in $X$ and $\ce$ a Dirichlet form on $L_2(\Omega)$.
Then the capacity of subsets $A\subset \overline\Omega$ and in particular subsets of the boundary $\Gamma={\overline\Omega}\backslash\Omega$ 
of $\Omega$ can be defined similarly.
One sets
\begin{eqnarray*}
\capp_\ce(A)&=&\inf\{\|\varphi\|_{D(\ce)}^2:\varphi\in D(\ce)  \mbox{ and there exists an open set  } \\[5pt]
&&\hspace{3cm}{}U\subset X \mbox{ with } A\subset U \mbox{ and } \varphi=1 \mbox{ on } U\cap \Omega\,\}
\;.
\end{eqnarray*}
Now the second criterion for uniqueness is formulated in terms of the capacity of the boundary of $\Omega$.

\begin{prop} \label{pcap3.2} Let $\Omega$ be a bounded,  open, connected subset of $X$
 with boundary  $\Gamma$.
Further let $\ce$ be a strongly local, inner regular,  Dirichlet form on $L_2(\Omega)$ and $\ce_m$, $\ce_M$ the corresponding extremal forms.
Then the  following conditions are equivalent:
\begin{tabel}
\item\label{pcap3.2-1}
$\ce_m=\ce_M$,
\item\label{pcap3.2-2}
$\capp_{\ce_m}(\Gamma)=0$.
\end{tabel}
\end{prop}
\proof\ \noindent\ref{pcap3.2-1}$\Rightarrow$\ref{pcap3.2-2}.$\;$
It follows from Corollary~\ref{ccap3.1} that there exists a sequence $\eta_n$  satisfying Condition~\ref{ccap3.1-2} of the corollary.
Since $\eta_n\in C_c(\Omega)$ it  follows that $\Gamma\subset (\supp\eta_n)^{\rm c}$.
Set $\varphi_n=\one_\Omega-\eta_n$.
Then $0\leq \varphi_n\leq 1$, $\varphi_n=1$ on $(\supp\eta_n)^{\rm c}$ and $\|\varphi_n\|_{D(\ce_m)}\to0$
as $n\to\infty$.
Therefore $\capp_{\ce_m}(\Gamma)=0$.

\smallskip

\noindent\ref{pcap3.2-2}$\Rightarrow$\ref{pcap3.2-1}.$\;$
Since $\capp_{\ce_m}(\Gamma)=0$ there exists a sequence of open sets $U_n\subset X$ such that $\Gamma\subset U_n\cap \,\overline\Omega$
and a sequence of $\varphi_n\in D(\ce_m)$ with $0\leq \varphi_n\leq 1$, $\varphi_n=1$ on $U_n\cap\Omega$ and $\|\varphi_n\|_{D(\ce_m)}\to0$ as
$n\to\infty$.
Now set $\eta_n=\one_\Omega-\varphi_n$.
Then $\one_\Omega-\eta_n\in B(\ce_m)$ and $\|\one_\Omega-\eta_n\|_{D(\ce_m)}\to0$ as $n\to\infty$.
But  $\supp\eta_n\subseteq  U_n^{\rm c}\cap \Omega$. 
Therefore $\eta_n\in B(\ce_M)$.
Thus the sequence of $\eta_n$ satisfies Condition~\ref{ccap3.1-2}$^\prime$ of Corollary~\ref{ccap3.1} which is equivalent to $\ce_m=\ce_M$.
\hfill$\Box$

\bigskip

The foregoing proof  is based on the assumption that  $\Omega$ is bounded.
But boundedness is  probably not essential for the conclusion of the proposition. 
The comparable result for the forms associated with elliptic differential operators on $\Ri^d$ is valid for unbounded domains
but the proof is considerably more complicated and depends on subadditivity of the capacity and localization arguments.

\end{document}